\documentclass[a4paper, 12pt]{article}

\usepackage{float}
\usepackage[sort&compress]{natbib}
\bibpunct{(}{)}{;}{a}{}{,} 
\usepackage{rotating}
\usepackage{lscape}
\usepackage{natbib}
\usepackage{bbm}
\usepackage{amsthm, amsmath, amssymb, mathrsfs, multirow, url,bm}
\usepackage{graphicx} 
\usepackage{ifthen} 
\usepackage{amsfonts}
\usepackage[usenames]{color}
\usepackage{fullpage}
\usepackage{caption}
\usepackage{subcaption}
\theoremstyle{plain} 
\newtheorem{thm}{Theorem}

\newtheorem{lem}{Lemma}

\theoremstyle{definition}

\newtheorem{model}{Model}

\theoremstyle{remark}
\newtheorem{remark}{Remark}

\newcommand{\calA}{\mathcal{A}}

\newcommand{\RR}{\mathbb{R}}

\newcommand{\E}{\mathsf{E}}

\newcommand{\Bs}{B^\star}
\newcommand{\eps}{\varepsilon}

\newcommand{\nm}{\mathsf{N}}
\newcommand{\N}{\mathsf{N}}
\newcommand{\bY}{{Y}}

\newcommand{\unif}{\mathsf{Unif}}

\newcommand{\btheta}{{\theta}}
\newcommand{\bbeta}{{\beta}}
\renewcommand{\S}{\mathcal{S}}

\newcommand{\B}{\mathcal{B}}
\newcommand{\BB}{\mathbb{B}}

\newcommand{\ZBs}{Z_{B^\star(s)}}

\newcommand{\one}{\mathbbm{1}}

\title{Empirical priors and posterior concentration in a piecewise polynomial sequence model}
\author{
Chang Liu\footnote{Department of Statistics, North Carolina State University} \qquad Ryan Martin$^*$ \qquad Weining Shen\footnote{Department of Statistics, University of California Irvine}
}
\date{\today}

\begin{document}
 
\maketitle 

\begin{abstract}
Inference on high-dimensional parameters in structured linear models is an important statistical problem.  This paper focuses on the case of a piecewise polynomial Gaussian sequence model, and we develop a new empirical Bayes solution that enjoys adaptive minimax posterior concentration rates and improved structure learning properties compared to existing methods. 
Moreover, thanks to the conjugate form of the empirical prior, posterior computations are fast and easy.  Numerical examples also highlight the  method's strong finite-sample performance compared to existing methods across a range of different scenarios.   

\smallskip

\emph{Keywords and phrases:} Bayesian estimation; change-point detection; high dimensional inference; structure learning; trend filtering.
\end{abstract}

\section{Introduction}
\label{S:intro}

Consider a Gaussian sequence model 
\begin{equation}
\label{eq:model}
Y_i \sim \nm(\theta_i, \sigma^2), \quad i=1,\ldots,n,
\end{equation}
where $Y=(Y_1,\ldots,Y_n)^\top$ are independent, the variance $\sigma^2 > 0$ is known, and inference on the unknown mean vector $\theta = (\theta_1,\ldots,\theta_n)^\top$ is desired.  It is common to assume that $\theta$ satisfies a {\em sparsity structure}, i.e., most $\theta_i$'s are zero, and work on these problems goes back at least to \citet{donoho1994minimax}, and more recently in  \citet{johnstone2004needles}, \citet{jiang2009general}, \citet{castillo2012needles},
{\citet{martin2014asymptotically}}, 
\citet{van2017adaptive},  {\citet{martin2020empirical}}, 
etc. 

{
There has also been recent interest in imposing different low-dimensional structures on high-dimensional parameters, namely, {\em piecewise constant} and, more generally, {\em piecewise polynomial}.  For a fixed positive integer $K$, we say that the $n$-vector $\theta$ has a piecewise degree-$K$ polynomial structure if there exists a simple partition $B$ of the index set into consecutive blocks $B(s) \subseteq \{1,\ldots,n\}$, with $s=1,\ldots,|B|$, such that, for each block $B(s)$, the corresponding sub-vector $\{\theta_j: j \in B(s)\}$ can be expressed as a degree-$K$ polynomial of the indices $j \in B(s)$.  This piecewise polynomial form is determined by the degree $K$ and the complexity $|B|$ of the block, i.e., its dimension is $(K+1)|B|$.  When this number is smaller than $n$, then a $\theta$ of this form clearly has a relatively low-dimensional structure.  For example, the piecewise constant case corresponds to $K=0$, so the complexity is completely determined by the number of blocks $|B|$.  
}

{Compared to sparse Gaussian signals, there is limited  literature studying piecewise constant and piecewise polynomial Gaussian sequence models. Regularization methods such as trend filtering \citep{kim2009ell_1} and locally adaptive regression splines \citep{mammen1997locally} are proposed to estimate the signal adaptively and recover the underlying block partitions. For piecewise constant problems, \citet{tibshirani2005sparsity} introduce  fused lasso based on a penalized least squares problem using the total variation penalty. \citet{rinaldo2009properties}, \citet{qian2016stepwise} investigate convergence rate of the fused lasso estimator and the asymptotic properties of pattern recovery. For signals with a more general piecewise polynomial structure,  \citet{tibshirani2014adaptive} propose adaptive piecewise polynomial estimation via trend filtering through minimizing a penalized least squares criterion, in which the penalty term sums the absolute $K$th order discrete derivatives over input points. \citet{guntuboyina2020adaptive} show that, under the strong sparsity setting and minimum length condition, the trend filtering estimator achieves  $n^{-1}$-rate, up
to a logarithmic multiplicative factor. In Bayesian domain,  methods such as Bayesian fused lasso \citep{kyung2010penalized} and Bayesian trend filtering \citep{roualdes2015bayesian} are proposed accordingly. However, to the best of our knowledge, no Bayesian literature has covered posterior contraction of adaptive estimation and asymptotic structure recovery for such piecewise polynomial Gaussian sequence models. Our goal here is to fill this gap.
}

{
Given the relatively low-dimensional representation of the high-dimensional $\theta$, the now-standard Bayesian approach would be to assign a prior for the unknown block configuration $B$ and a conditional prior on the block-specific $(K+1)$-dimensional parameters that determine the polynomial form.  For the prior on $B$, the goal would be to induce ``sparsity'' in the sense that the prior concentrates on block configurations $B$ with $|B|$ relatively small.  For this, one can mostly follow the existing Bayesian literature on sparsity structures, e.g., \citet{castillo2012needles}, \citet{castillo2015bayesian}, \citet{martin2017empirical}, \citet{liu2018bayesian}, etc.  However, for the quantities that determine the polynomial form on a given block configuration, the situation is very different.  In the classical sparsity settings, it is reasonable to assume that those signals that are not exactly zero are still relatively small, so a conditional prior centered around zero can be effective.  In this piecewise polynomial setting, there is no obvious fixed center around which a prior should be concentrated.  Of course, one option is to choose a fixed center and wide spread, but then the tails of the prior distribution become particularly relevant (e.g., \citet{castillo2012needles}, Theorem~2.8).  An alternative is to follow \citet{martin2019data}, building on \citet{martin2014asymptotically} and \citet{martin2017empirical}, and use an {\em empirical prior} that lets the data help with correctly centering the prior distribution.  
}

Details of this empirical prior construction are presented in Section~\ref{S:model}.  
Our theoretical results in Section~\ref{S:rate} demonstrate that the corresponding empirical Bayes posterior distribution enjoys adaptive concentration at the same rate of trend filtering, adjusting to phase transitions, but requires weaker conditions than that in \citet{guntuboyina2020adaptive}. In addition, we establish structure learning consistency results which, to our knowledge, is the first one for piecewise polynomial sequence models in the Bayesian literature.
  Furthermore, since the proposed empirical priors are conjugate, the posterior is relatively easy to compute, and the numerical simulations in Section~\ref{S:examples} which compares our method with trend filtering,  demonstrate the advantageous performance of our method in signal estimation and structure recovery under finite-sample settings. In Section~\ref{S:real}, we apply our method to two real-world applications where the underlying truths are considered to be piecewise constant and piecewise linear respectively. Finally, some concluding remarks are made in Section~\ref{S:discuss}, and technical details and proofs are presented in the Appendix.

\section{Empirical Bayes formulation}
\label{S:model}

\subsection{Piecewise polynomial model}
\label{SS:formulation}

Before we can introduce our proposed prior and corresponding empirical Bayes model, we need to be more precise about the within-block polynomial formulation.  Start with the case $|B|=1$ corresponding to there being only one block.  A vector $\theta$ being a degree-$K$ polynomial with respect to $B$ corresponds to $\theta \in \S$, where 
\begin{equation}
\label{eq:span}
\S = \text{span}\{v_0,v_1,\ldots,v_K\}, 
\end{equation}
and $v_k = (1^k, 2^k,\ldots, n^k)^\top \in \RR^n$, with $k=0,1,\ldots,K$.  In other words, if $Z \in \RR^{n \times (K+1)}$ is a matrix whose columns form a basis for $\S$, then $\theta$ can be expressed as $Z \beta$ for some vector $\beta \in \RR^{K+1}$.  More generally, for a generic simple partition $B$, if $\theta$ is a piecewise degree-$K$ polynomial on the block configuration $B$ as described in Section~\ref{S:intro}, then it can be expressed as $Z^B \beta^B$, where 
\begin{equation}
    Z^B=
\begin{pmatrix}
Z_{B(1)} & 0 &\ldots & 0\\
0 & Z_{B(2)} & \ldots & 0\\
\vdots & \vdots &\ldots &\vdots\\
0 &0 &\ldots&Z_{B(|B|)}
\end{pmatrix} \in \RR^{n\times |B|(K+1)}, 
\label{eq:Z^B}
\end{equation}
$Z_{B(s)}$ is the sub-matrix of $Z$ with its row indices included in $B(s)$, and
\begin{equation}
    \beta^B = \begin{pmatrix}
\beta_1^B\\
\vdots\\
\beta_{|B|}^B
\end{pmatrix} \in \RR^{|B|(K+1)},
\quad \beta_s^B \in \RR^{K+1}, \quad s=1,\ldots,|B|.
\label{eq:theta^B}
\end{equation}
The following two examples will illustrate the piecewise polynomial formulation.
\begin{itemize}
\item When $K=0$, the vector $\theta$ formed by $Z^B\beta^B$ is piecewise constant. For a specific block segment $B(s)$, we can write $Z_{B(s)}=\one_{|B(s)|}$ and, therefore, \[\theta_i \equiv \beta^B_s \in \RR, \quad i \in B(s).\]
Note that in this case the Gaussian sequence model can be rewritten in the form of a one-way analysis of variance model with $|B|$ treatments and $|B(s)|$ number of replications in each treatment, $s=1,\ldots,|B|$.
\vspace{-2mm}
\item When $K=1$, the vector $\theta$ formed by $Z^B\beta^B$ is piecewise linear. For a specific block segment $B(s)$, we can write
\[Y_{B(s)}=Z_{B(s)}\beta^B_s+\eps, \quad \eps\sim \N_{|B(s)|}(0, \sigma^2 I),\]
where $Y_{B(s)}$ is the sub-vectors of $Y$ with its indices in $B(s)$, and $\beta^B_s$ is a two-dimensional vector. Hence, within each segment, the observed data can be viewed as a random sample generated from a block-specific simple linear regression model with intercept and slope being $\beta^B_{s, 1}$ and $\beta^B_{s,2}$.
\end{itemize}



To summarize, if $\theta$ is a $n$-vector that is assumed to have a piecewise degree-$K$ polynomial structure, then we can reparametrize $\theta$ as $(B,\theta^B)$, where $\theta^B$ is expressed as $Z^B \beta^B$, for some $\beta^B \in \RR^{|B|(K+1)}$, and $Z^B$ is as in \eqref{eq:Z^B} for some generator matrix $Z$ whose columns form a basis for $\S$ in \eqref{eq:span}.  The matrix $Z$ is not unique and, therefore, $\beta^B$ is not unique either.  But interest is in the projection $\theta^B$, which is independent of the choice of basis, so this non-uniqueness will not be a problem in what follows.  


\subsection{Empirical prior}
\label{SS:prior}

In light of our representation of a piecewise polynomial mean vector $\theta$ via $(B,\theta^B)$ or $(B,\beta^B)$, a hierarchical representation of the prior distribution will be most convenient.  That is, we first specify a prior for $B$, then a conditional prior for $\beta^B$, given $B$; this in turn will induce a conditional prior for $\theta^B$.  Here we follow this general prior specification strategy, but with a slight twist wherein we allow the conditional prior for $\beta^B$ to depend on data in a particular way.  Then this empirical prior for $(B,\beta^B)$ will immediately induce a corresponding empirical prior for $(B,\theta^B)$ and, finally, for $\theta$.  

Intuitively, there is no reason to introduce a piecewise polynomial structure if not for a belief that there are not too many  blocks, i.e., that $|B|$ is relatively small compared to $n$; see Section~\ref{S:rate}.  This belief can be incorporated into the prior for $B$ in the following way.  Set $b=|B|$, and introduce a marginal prior  
\begin{equation}
\label{eq:S.prior}
f_n(b) \propto n^{-\lambda (b-1)}, \quad b=1,\ldots,n, 
\end{equation}
where $\lambda > 0$ is a constant to be specified.  Note that this is effectively a truncated geometric distribution with parameter $p=n^{-\lambda}$, which puts most of its mass on small values of the block configuration size, hence incorporating the prior information that $\theta$ is not too complex.  Next, if the configuration size $b$ is given, the blocks correspond to a simple partition of $\{1,2,\ldots,n\}$ into $b$ consecutive chunks, and there are $\binom{n-1}{b-1}$ such partitions.  So, for the conditional prior distribution of $B$, given $|B|$, we can take a discrete uniform distribution.  Therefore, the prior distribution for $B$ is given by
\begin{equation}
\label{eq:B.prior}
\pi_n(B) = \textstyle\binom{n-1}{|B|-1}^{-1} f_n(|B|), 
\end{equation}
where $B$ ranges over all simple partitions of $\{1,2,\ldots,n\}$ into consecutive blocks.  

{
Next we give the conditional prior for $\beta^B$, given $B$.  We propose to assign independent normal priors to each $\beta^B_s$ corresponding to segment $B(s)$.  A notable departure from the traditional Bayesian formulation is that we follow \citet{martin2017empirical} and let the data inform the prior center.  Specifically, the conditional prior for $\beta^B$, given $B$, is taken to be 
\begin{equation}
     {\beta}^B_{s} \sim \N_{K+1}(\hat{ {\beta}}^B_{s}, v\{ Z_{B(s)}^\top Z_{B(s)}\}^{-1}), \quad s=1,\ldots,|B|, \quad \text{independent},
    \label{eq:prior_beta}
\end{equation}
where $\hat{\beta}^B_s$ is the least-squares estimator 
\[ \hat\beta_s^B = (Z_{B(s)}^\top Z_{B(s)})^{-1} Z_{B(s)}^\top Y_{B(s)}, \]
and $v>0$ is a constant controlling prior spread.  Write the conditional density function of $\beta^B$, given $B$, with respect to Lebesgue measure on $\RR^{|B|(K+1)}$, as
\[ \tilde \pi_n(\beta^B \mid B) = \prod_{s=1}^{|B|} \nm_{K+1}(\beta_s^B \mid \hat\beta_s^B, v\{ Z_{B(s)}^\top Z_{B(s)}\}^{-1}), \]
a product of individual $(K+1)$-variate normal densities.  This induces a prior on $\theta^B$ through the mapping $\theta^B = Z^B \beta^B$ that defines it, but since this is generally not a bijection, there is no density function with respect to Lebesgue measure on $\RR^n$.  To see this, let $\theta_{B(s)}^B$ denote the sub-vector of $\theta^B$ with indices included in $B(s)$, then we can observe that the induced conditional prior on $\theta_{B(s)}^B$ is $\nm_{|B(s)|}(P_{B(s)} Y_{B(s)}, v P_{B(s)})$, where 
\begin{equation}
\label{eq:projection}
P_{B(s)} = Z_{B(s)}\{ Z_{B(s)}^\top Z_{B(s)}\}^{-1} Z_{B(s)}^\top
\end{equation}
is the matrix that projects onto the space spanned by the columns of $Z_{B(s)}$.  Since $P_{B(s)}$ is a projection, it is not full rank and, therefore, the prior for  $\theta_{B(s)}^B$ is a degenerate normal.  Despite this degeneracy, the conditional prior for $\theta^B$, given $B$, still exists; it is just more convenient to express in terms of the conditional prior for $\beta^B$.  That is, we define the conditional empirical prior for $\theta^B$, given $B$, as 
\[ \Pi_n(\calA \mid B) = \int_{\{\beta^B: Z^B \beta^B \in \calA\}} \tilde \pi_n(\beta^B \mid B) \, d\beta^B, \quad \calA \subseteq \RR^n. \]
Note that while the prior for $\beta^B$ depends on the particular basis in $Z$, the prior for $\theta^B$ only depends on the projection which does not depend on the choice of basis.  Finally, our empirical prior for $\theta$ is defined as 
\begin{align*}
\Pi_n(\calA) & = \sum_B \pi_n(B) \, \Pi_n(\calA \mid B) \\
& = \sum_B \pi_n(B) \int_{\{\beta^B: Z^B\beta^B \in \calA\}} \tilde \pi_n(\beta^B \mid B) \, d\beta^B.
\end{align*}
}

The reader may be anticipating that the combination of an empirical prior with the likelihood amounts to double-use of data.  To avoid potentially over-fitting, we propose the following mild additional regularization.  Let $L_n(\theta)$ denote the likelihood function based on the model \eqref{eq:model}, i.e., $L_n(\theta) \propto \exp\{-\tfrac{1}{2\sigma^2} \|Y-\theta\|^2\}$, where $\|\cdot\|$ denotes the $\ell_2$-norm on $\RR^n$.  For a fixed $\alpha \in (0,1)$, define a regularized empirical prior 
\begin{equation}
\label{eq:post00}
 \Pi_n^\text{reg}(d\theta) \propto L_n(\theta)^{-(1-\alpha)} \, \Pi_n(d\theta).
\end{equation}
Dividing by a fractional power of the likelihood effectively down-weights those parameter values with especially large likelihood, hence discouraging over-fitting.  Typically, one would take $\alpha$ to be close to 1---for example, we take $\alpha=0.99$ in the simulation examples presented in Section~\ref{S:examples}---so this additional regularization is very mild indeed.

\subsection{Posterior}
\label{SS:posterior}

For the posterior distribution, we propose to combine the regularized empirical prior $\Pi_n^\text{reg}$ with the likelihood $L_n$ according to Bayes's formula:
\begin{align}
    \label{eq:post0}
\Pi^n(\calA) &\propto \int_{\calA} L_n(\theta) \,  \Pi_n^\text{reg}(\theta)
\end{align}
The following sections investigate the theoretical convergence properties and practical performance of this empirical Bayes posterior distribution.  

Of course, the posterior in \eqref{eq:post0} can be rewritten as 
\[ \Pi^n(\calA) \propto \int_{\calA} L_n(\theta)^\alpha \, \Pi_n(d\theta), \]
which is particularly well-suited for our theoretical analysis.  This sort of {\em generalized Bayes posterior} has received considerable attention recently, e.g., \citet{grunwald2017inconsistency}, \citet{miller2019robust}, \citet{holmes2017assigning}, \citet{syring2019calibrating}, and \citet{bhattacharya2019bayesian}, though not specifically for the purpose of regularization.  One might ask if $\alpha=1$ is a valid choice, since this makes the above display look more like the familiar Bayesian update, but the answer is unclear because our analysis here makes specific use of $\alpha < 1$. The reason we work with $\alpha < 1$ is for simplicity, but there is nothing to gain by including $\alpha=1$.  In particular, we improve upon the existing Bayesian rate results for this problem (Remark~\ref{re:pas.rockova}), in some cases achieving optimal rates, and give new results on Bayesian structure learning. 
And even if the reader is uncomfortable with giving up a tiny portion of the likelihood, he/she can interpret $\Pi^n$ as the full likelihood combined with a data-dependent prior as in \eqref{eq:post0}, closer to traditional empirical Bayes.  

A practical benefit to the relative simplicity of our formulation is that the posterior distribution turns out to be not so complicated.  Indeed, by combining \eqref{eq:prior_beta} and \eqref{eq:post0}, the posterior distribution $\Pi^n$ for $\theta$ is given by 
\begin{equation}
\label{eq:post_theta}
\Pi^n(\calA) 
= \sum_{B} \pi^n(B) \int_{\{\beta^B:Z^B\beta^B \in \calA\}} \prod_{s=1}^{|B|} \N_{K+1}(\beta_s^B \mid \hat{\beta}_s^B, \tfrac{\sigma^2 v}{\sigma^2 + \alpha v}\{ Z_{B(s)}^\top Z_{B(s)}\}^{-1}) \, d\beta^B,
\end{equation}
where 
\begin{equation}
\pi^n(B)\propto \pi_n(B)(1+\tfrac{v\alpha}{\sigma^2})^{-(K+1)|B|/2}e^{-\frac{\alpha}{2\sigma^2}\sum_{s=1}^{|B|}\|(I-P_{B(s)})Y_{B(s)}\|^2},
\label{eq:marginal_B}
\end{equation}
with $P_{B(s)}$ the projection in \eqref{eq:projection}.  From this expression, we can see that there are three major factors contributing to the log-marginal posterior distribution of $B$: the prior distribution for block configuration $\log \pi_n(B)$, a penalty term on model complexity proportional to $-|B|$, and a model fitting measure proportional to the negative sum of squared residuals.  Therefore, our posterior distribution would prefer models with fewer blocks and better fitting given the observed data $Y$. Details about how we compute the posterior distribution are presented in Section~\ref{S:computation}.  


\section{Asymptotic properties}
\label{S:rate}

\subsection{Setup}
\label{SS:setup}

For a vector $\theta \in \RR^n$ that has a piecewise degree-$K$ polynomial structure, write $B_\theta$ for its block configuration, and let $|B_\theta|$ denote its cardinality.  Then our parameter space corresponds to $\Theta_n(K)$, the set of all $n$-vectors with a piecewise degree-$K$ polynomial structure and having $|B_\theta| = o(n)$.  The latter condition on the size of the block configuration ensures that there are not too many blocks, i.e., that the signal is not too complex.  

\begin{figure}[t]
\centering
    \begin{subfigure}[b]{0.5\textwidth}
                \centering
                \includegraphics[width=0.85\textwidth]{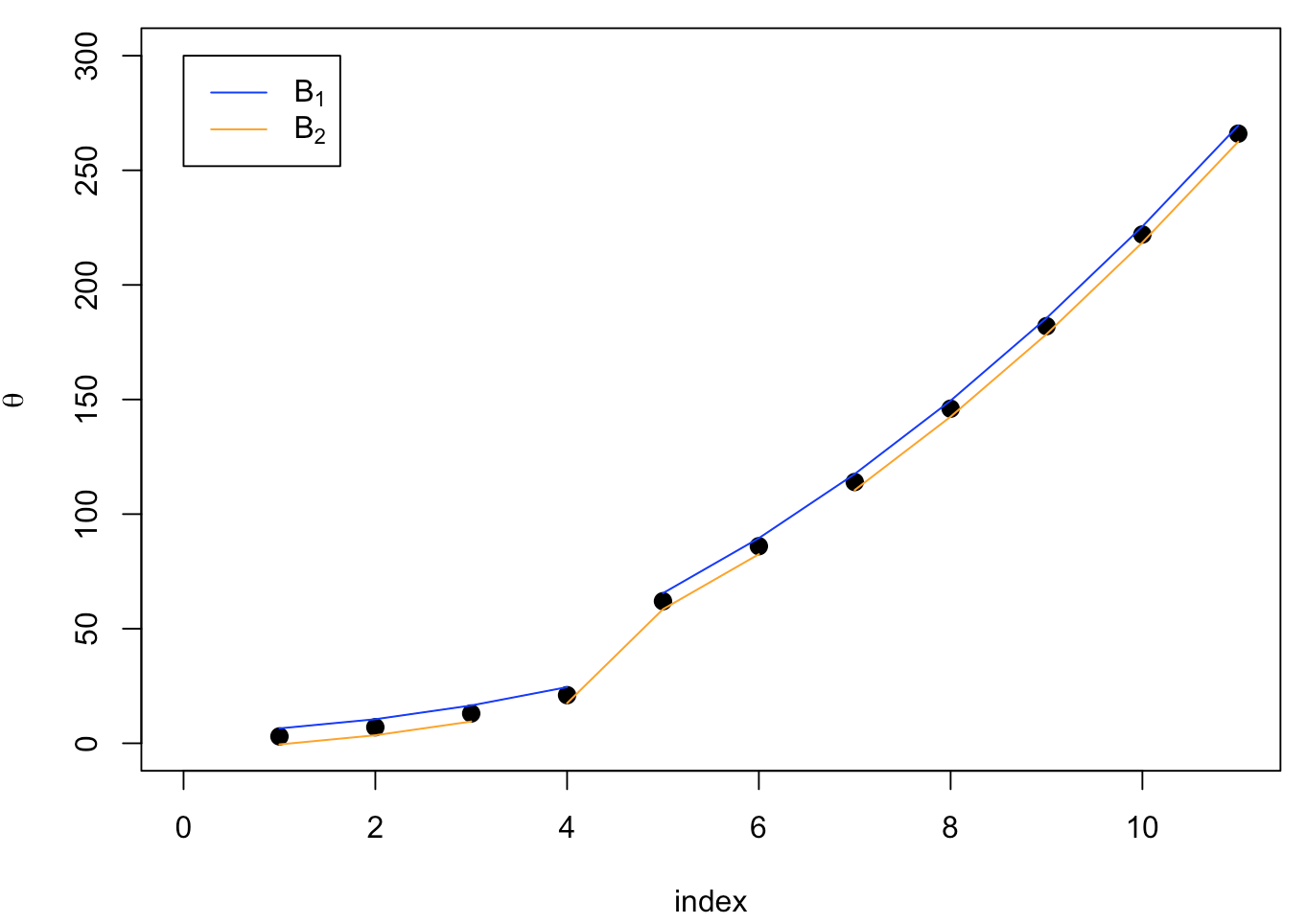}
                \caption{Can be resolved}
    \end{subfigure}%
    \begin{subfigure}[b]{0.5\textwidth}
                    \centering
                \includegraphics[width=0.85\textwidth]{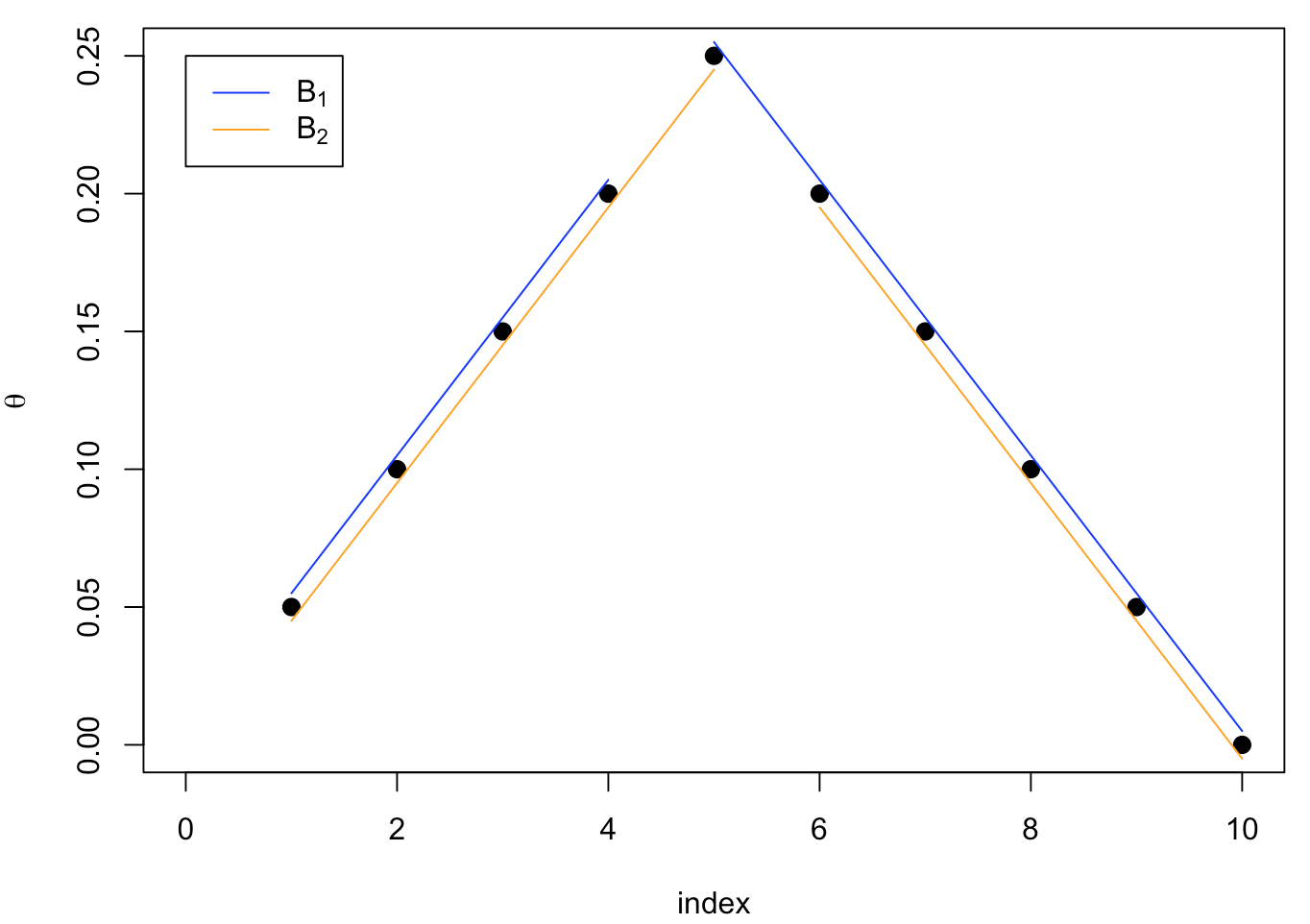}
                \caption{Cannot be resolved}
    \end{subfigure}%
\caption{Two examples pertaining to  identifiability of $B_\theta$ for a given  $\theta$, one where the non-uniqueness of $B_\theta$ can be resolved and one where it cannot. The two lines on each plot both pass through the marked points; the small jitter is to help distinguish the groupings corresponding to the two block configurations.}
\label{fig:identifiablity}
\end{figure}

When $K \geq 1$, it is possible that a vector $\theta$ has multiple block configurations $B_\theta$.  That is, there could be multiple $B$ and $\beta^B$ such that $\theta = Z^B \beta^B$.  This does not present a problem for questions related to estimation of $\theta$, but it does create some identifiability concerns in the context of structure learning, i.e., recovering the underlying block structure.  In some cases, the non-uniqueness can be resolved by defining $B_\theta$ as the ``most economical'' of the candidate $B$'s.  For example, for an arbitrary signal vector, any $(K+1)$-tuple of consecutive points could be fit perfectly by a degree-$K$ polynomial, so blocks of size $K+1$ or smaller are meaningless and should be ruled out.  Panel~(a) of Figure~\ref{fig:identifiablity} shows an illustration of this for $K=2$, piecewise quadratic.  However, there are other cases where the non-uniqueness cannot be resolved by ruling out blocks that are too small.  Panel~(b) of Figure~\ref{fig:identifiablity} shows an example of this, where the two candidate block configurations are perfectly indistinguishable by data.  Again, this is of no concern for results in Section~\ref{SS:recovery} below, so we postpone our discussion of how this is resolved until Section~\ref{SS:structure}. 


\subsection{Posterior concentration rates}
\label{SS:recovery}

For $x \in \RR^n$, define the scaled $\ell_2$-norm $\|x\|_n=n^{-1/2}\|x\|_2$ and, for $\theta^\star \in \Theta_n(K)$, define  
\begin{equation}
\label{eq:target}
\eps_n^2(\theta^\star) = \begin{cases} n^{-1} & \text{if $|B_{\theta^\star}|=1$}, \\
n^{-1} |B_{\theta^\star}| \log n & \text{if $|B_{\theta^\star}| \ge 2$}. \end{cases}
\end{equation}
Note that, in the case $|B_{\theta^\star}|=1$, the best estimator of $\theta^\star$ would be $P_{\S} Y$, where $P_{\S}$ is the projection matrix onto $\S$ in \eqref{eq:span}, and its expected sum-of-squared-error is $O(n^{-1})$, consistent with \eqref{eq:target}.  For the case with $|B_{\theta^\star}|\ge 2$,  the rate \eqref{eq:target} is consistent with others obtained in the literature; see Remark~\ref{re:guntuboyina} below.  
Theorem~\ref{thm:rate} says that the $\Pi^n$ constructed above attains the rate defined in \eqref{eq:target}.  Since the prior can achieve the rates $\eps_n^2(\theta^\star)$ defined above without knowledge of $\theta^\star$ or $|B_{\theta^\star}|$, it follows that our posterior concentration results are {\em adaptive} to the unknown complexity of $\theta^\star$.   




\begin{thm}
\label{thm:rate}
Consider the model \eqref{eq:model} with known $\sigma^2 > 0$ and assume that $\theta^\star$ has a piecewise polynomial structure of degree $K \geq 0$, with $K$ known.  Let $\Pi^n$ be the corresponding empirical Bayes posterior distribution for $\theta \in \RR^n$ described above.  If $\eps_n^2(\theta^\star)$ is as in \eqref{eq:target}, then for any sequence $M_n$ with $M_n \to \infty$, there exists a constant $G > 0$ such that 
\[ \E_{\theta^\star} \Pi^n(\{\theta \in \RR^n: \|\theta-\theta^\star\|^2_n > M_n \eps_n^2(\theta^\star)\}) \lesssim e^{-G M_n n \eps_n^2(\theta^\star)}. \]
for all large $n$, uniformly in $\theta^\star \in \Theta_n(K)$.  For the latter case in \eqref{eq:target}, the sequence $M_n$ above can be replaced by a sufficiently large constant $M > 0$.  
\end{thm}

\begin{remark}
\label{re:oracle}
Given data $Y \sim \N_n(\theta^\star
,\sigma^2 I)$, an oracle who has access to $B_{\theta^\star}$ would fit a polynomial of degree $K$ in each of the partitions given by $B_{\theta^\star}$. This would be a linear estimator and its corresponding oracle risk is $O(n^{-1}|B_{\theta^\star}|)$.  Note that the rate achieved in Theorem~\ref{thm:rate} is comparable to the oracle risk.  Indeed, our method adaptively learns the underlying block structure of $\theta^\star$ and, in the case $|B_{\theta^\star}|=1$, we can exactly match the oracle rate; otherwise, the price we pay in terms of the rate is only logarithmic.  
\end{remark}

\begin{remark}
\label{re:guntuboyina}
The minimax optimal rate, $n^{-1} |B_{\theta^\star}| \log(en/|B_{\theta^\star}|)$, can be achieved if more control on the complexity of $\theta^\star$ is assumed, namely, if $|B_{\theta^\star}| = O(n^t)$ for some $t \in [0,1)$.  The only way this extra assumption fails is if the signal is extremely complex, e.g., if $|B_{\theta^\star}|=O(n/\log n)$.  Such cases effectively have no low-dimensional block structure and should be rare in practice.  This minimax rate can be achieved by trend filtering \citep[see][Corollary~2.3]{guntuboyina2020adaptive}, but this too requires additional assumptions.  
Indeed, their result holds only when their minimum length condition is satisfied and the tuning parameter is properly chosen within an unspecified ``ideal'' range. The former---see Equation~(13) in \citet{guntuboyina2020adaptive}---restricts the length of the minimal block to be no smaller than $O(n |B_{\theta^\star}|^{-1})$, which cannot be checked in applications.  They also make a strong sparsity assumption that requires $|B_{\theta^\star}|$ to be ``much smaller than $n$.'' This surely excludes extremely high-complexity cases like $|B_{\theta^\star}|=O(n/\log n)$. Therefore, our empirical Bayes posterior concentration rate result is no weaker than the results for trend filtering in \citet{guntuboyina2020adaptive} which those authors argue are stronger than any existing results in the literature. 
\end{remark}

\begin{remark}
\label{re:pas.rockova}
\citet{van2017bayesian} present a result similar to that in Theorem~\ref{thm:rate}, for the piecewise constant case $K=0$, with a rate of $|B_{\theta^\star}| \log(n / |B_{\theta^\star}|)$.  However, translating their notation to ours, they assume bounds on both $\|\theta^\star\|$ and on $|B_{\theta^\star}|$, which we do not require.  And in light of Theorem~2.8 of \citet{castillo2012needles}, we do not expect that optimal concentration rates can be achieved using their fixed-center normal prior for $\theta^B$, given $B$, without some assumptions on the magnitude of $\theta^\star$.  
\end{remark}


Next, we show that the posterior mean $\hat\theta = \int \theta \, \Pi^n(d\theta)$ is an adaptive, asymptotically minimax estimator.

\begin{thm}
\label{thm:mean}
Under the setup in Theorem~\ref{thm:rate}, 
\[ \E_{\theta^\star} \|\hat\theta - \theta^\star\|_n^2 \lesssim M_n \eps_n^2(\theta^\star), \]
for all large $n$, uniformly in $\theta^\star \in \Theta_n(K)$.  In the latter two cases of \eqref{eq:target}, the diverging sequence $M_n$ can be replaced by a constant $M$ which can be absorbed into ``$\lesssim$'' above. 
\end{thm}

\subsection{Structure learning}
\label{SS:structure}

In addition to estimation consistency, it is interesting to consider when the posterior is able to recovery the block structure of the true piecewise polynomial signal $\theta^\star$.  To our knowledge, this is the first Bayesian (or empirical Bayesian) investigation into structure learning in the piecewise polynomial Gaussian sequence model. When $K=0$, i.e., the true signal is piecewise constant, learning the underlying block structure can be interpreted as detection of the ``change points" or ``jump points", which has many real-world applications. In the non-Bayesian literature, structure recovery for piecewise constant and piecewise polynomial signals has received some  attention, and below we compare our results with those available for the trend filtering, binary segmentation, etc. 

As a first result in this direction, Theorem~\ref{thm:dim} says that the effective dimension of the posterior is no larger than a multiple of the true block configuration size---in other words, the posterior is of roughly the correct complexity.  Note that this result only pertains to the size $|B_\theta|$ of the block configurations, which can be uniquely determined, so there are no identifiability issues here.  Finally, for this and the other results of this section, the statements are formulated in in terms of the marginal posterior distribution $\pi^n$ for the block configuration $B$, as defined in \eqref{eq:marginal_B}.  

\begin{thm}
\label{thm:dim}
Under the setup in Theorem~\ref{thm:rate}, for any $C > 1 + \lambda^{-1}$, where $\lambda$ is as in \eqref{eq:S.prior}, there exists a constant $G > 0$ such that 
\[ \E_{\theta^\star} \pi^n(\{B: |B| > C |B_{\theta^\star}|\}) \lesssim e^{-G |B_{\theta^\star}|\log n}, \]
for all large $n$, uniformly in $\theta^\star \in \Theta_n(K)$.  
\end{thm}

Block configuration size is important, but we also care about identifying the underlying block structure.  Of course, before we can say any more about this, we need to address the potential non-identifiability of $B_{\theta^\star}$.  As we mentioned before, there are no such issues in the piecewise constant case with $K=0$, but non-identifiability is possible for other $K \geq 1$ cases.  On the one hand, if $\theta^\star$ is such that non-uniqueness can be resolved simply by taking the most economical of those equally-well-fitting block configurations, then that is how $B_{\theta^\star}$ is defined.  On the other hand, if $\theta^\star$ has multiple block configurations of the same size, like in Figure~\ref{fig:identifiablity}(b), then it is impossible to distinguish between these.  In such cases, the best we can hope for is that the posterior distribution will concentrate on the set $\BB^\star = \{B_{\theta^\star}\}$ of equivalent block configurations corresponding to $\theta^\star$ and, in fact, this is what the results below establish.   

The first result below concerns the event that $B$ is a {\em refinement} of $B_{\theta^\star}$, denoted by $B \sqsupset B_{\theta^\star}$, for some $B_{\theta^\star} \in \BB^\star$.  That is, if $B \sqsupset B_{\theta^\star}$, then every block in $B_{\theta^\star}$ can be expressed as a union of blocks in $B$ or, equivalently, no block in $B$ intersects with more than one block in $B^\star$.  Since refinements, or unnecessary splits of $B_{\theta^\star}$, are a sign of inefficiency, we hope that the posterior will discourage such cases.  Indeed, Theorem~\ref{thm:no.refinements} below shows that the posterior distribution assigns vanishing probability to the event ``$B \sqsupset B_{\theta^\star}$'', which means that the posterior for $B$ asymptotically avoids those inefficient refinements.  This is analogous to the ``no supersets'' theorems in \citet[][Thereom~4]{castillo2015bayesian} and \citet[][Theorem~4]{martin2017empirical} for variable selection in linear regression context.  The only additional requirement here is that the power $\lambda$ in the prior for $|B|$ in \eqref{eq:S.prior} not be too small; otherwise, the prior does not sufficiently penalize those block configurations that are too complex, leaving open the possibility for over-fitting.  Similar conditions appear in the regression setting, e.g., the conditions of Theorem~4 in \citet{castillo2015bayesian}.    
\begin{thm}
\label{thm:no.refinements}
Under the setup of Theorem~\ref{thm:rate}, 
\[ \E_{\theta^\star} \pi^n(\{B: B \sqsupset B_{\theta^\star} \text{ for some $B_{\theta^\star} \in \BB^\star$}\}) \to 0, \quad n \to \infty,  \]
uniformly in $\theta^\star \in \Theta_n(K) \cap \{\theta: |B_{\theta}| = o(n^\lambda)\}$, with $\lambda > 0$ as in the prior \eqref{eq:S.prior}.
\end{thm}

If $\lambda \geq 1$, then the above  condition on $|B_{\theta^\star}|$ is satisfied for all $\theta^\star \in \Theta_n(K)$.  However, for smaller values of $\lambda$, like those having good empirical performance in Section~\ref{S:examples}, restricting to a proper subset of the parameter space is required, but this is not severe.  

A natural follow-up question is if the true block configuration $B_{\theta^\star}$ or, more generally, the set $\BB^\star$ of equivalent true block configurations can be recovered exactly.  Before stating our affirmative answer to this question, we need some additional notation.  First, define the $0^\text{th}$- and $1^\text{st}$-order difference operators as $\Delta^0x = x$ and 
\[ \Delta^1 x = (x_2-x_1,x_3-x_2,\ldots,x_n-x_{n-1})^\top, \]
respectively, where $x \in \RR^n$.  For a generic order $K \geq 2$, the $K^\text{th}$-order difference, $\Delta^K: \RR^n \to \RR^{n-K}$, is defined recursively as $\Delta^K x = \Delta^1(\Delta^{K-1} x)$.  Second, a change in the signal $\theta^\star$ from one block to another can only be detected if the change is sufficiently large, and the definitions of ``change'' and ``sufficiently large'' are related to properties of the difference operators applied to $\theta^\star$.  In particular, the set of indices where a change in the $(K+1)^\text{st}$-order occurs is defined as 
\[ J_{\theta^\star} = \{j=1,\ldots,n-K-1: (\Delta^{K+1}\theta^\star)_j \neq 0\}. \]
In the piecwise constant case, with $K=0$, the set $\{j+1: j \in J_{\theta^\star}\}$ consists of those indices at which the signal jumps from one value to another.  Then both the minimal change in $\theta^\star$ on $J_{\theta^\star}$ and the minimal spacing between changes will be relevant to determining whether a change is sufficiently large to be detectable.  These are defined, respectively, as 
\[ \delta_n(\theta^\star) = \min_{j \in J_{\theta^\star}} |(\Delta^{K+1} \theta^\star)_j| \quad \text{and} \quad \gamma_n(\theta^\star) = \min_{j,j' \in J_{\theta^\star}, j \neq j'} |j-j'|. \]
Then the following theorem states that the block configuration $B_{\theta^\star}$ can be recovered exactly if $\gamma_n(\theta^\star) \delta_n^2(\theta^\star)$ is sufficiently large, analogous to the so-called {\em beta-min} condition in linear regression \citep[e.g.,][Chapter~2]{buhlmann2011statistics}.

{
\begin{thm}
\label{thm:B_consist}
Under the setup in Theorem~\ref{thm:rate}, suppose that 
\begin{equation}
    \label{eq:select}
    \gamma_n(\theta^\star)\delta^2_n(\theta^\star)\ge \frac{4^{K+1}M\sigma^2}{\alpha(1-\alpha)}\log n
\end{equation}
with $M>4+\lambda$ and $\lambda \ge 3$, where $\lambda$ controls the prior \eqref{eq:B.prior}.  Then 
\begin{equation}
\label{eq:B_consist}
\E_{\theta^\star} \pi^n(\BB^\star) \to 1, \quad n \to \infty. 
\end{equation}
\end{thm}
}

To our knowledge, only the piecewise constant ($K=0$) case---where the true $B_{\theta^\star}$ is unique---has been considered in the literature, so we focus on that version here in our discussion of Theorem~\ref{thm:B_consist}.  In that case, $\gamma_n(\theta^\star)$ and $\delta_n(\theta^\star)$ represent the smallest number of indices between jumps and the smallest signal jump in $\theta^\star$.  To make a parallel between the piecewise constant signal problem and a one-way analysis of variance, $\gamma_n(\theta^\star)$ is like the minimum number of replications across all the treatment groups and $\delta_n(\theta^\star)$ is like the minimum effect size.  In that classical analysis of variance context, where the number of treatment groups and group memberships are fixed and known, the F-test has power converging to 1 if $\gamma_n(\theta^\star) \delta_n^2(\theta^\star)$ is bounded away from 0.  The condition $\gamma_n(\theta^\star) \delta_n^2(\theta^\star) \gtrsim \log n$ in \eqref{eq:select} here is only slightly stronger, i.e., we only pay a logarithmic price for not knowing the number of groups or group memberships.  Returning to the general piecewise constant case, if the minimum block size $\gamma_n(\theta^\star)$ is fixed as $n$ and $|B_{\theta^\star}|$ go to infinity, the result in Theorem~\ref{thm:B_consist} above matches the pattern recovery property for fused lasso in \citet{qian2016stepwise}, and is stronger than the corresponding results in \citet{lin2017sharp} and \citet{dalalyan2017prediction}.  We can also allow the minimum block size $\gamma_n(\theta^\star)$ to grow. For example, the minimum block length condition in \citet{guntuboyina2020adaptive} states that $\gamma_n(\theta^\star)$ can be of order $O(n|B_{\theta^\star}|^{-1})$, corresponding to equally partitioning over blocks. In this case, the minimum jump size simply needs to satisfy,
 \[\delta^2(\theta^\star) \gtrsim n^{-1}|B_{\theta^\star}|\log n,\] 
which is mild since the right-hand side would typically be vanishing.  This flexibility makes our result preferable to those for fused lasso and comparable to that for wild binary segmentation in Theorem~3.2 of \citet{fryzlewicz2014wild}, which is the best result available in the literature that we are aware of.  Finally, we want to emphasize, again, that Theorems~\ref{thm:no.refinements}--\ref{thm:B_consist} are, to our knowledge, the first such results in the Bayesian literature.

\section{Computation}
\label{S:computation}

Genuine Bayesian solutions to high-dimensional problems, ones for which optimal posterior rates are available, tend to be based on non-conjugate, heavy-tailed priors, making computation non-trivial.  Our empirical Bayes solution, on the other hand, is based on a conjugate prior for $\theta^B$, making computations relatively simple.  


Indeed, recall that the marginal posterior for $B$ is available in closed-form, up to proportionality, as in \eqref{eq:marginal_B}.  Furthermore, recall from \eqref{eq:post_theta} that the conditional distribution of $\theta^B$, given $B$, is determined by a linear transformation of a normal random variable, which is easy to simulate.  Together, these two observations suggest the following Metropolis--Hastings algorithm to draw Markov chain Monte Carlo (MCMC) samples from the proposed posterior $\Pi^n$ for $\theta$:
\begin{enumerate}
    \item At iteration $t$, given current block partition $B^{(t)}$, sample $B' \sim q(\cdot \mid B^{(t)})$.
    \item Sample $U \sim \unif(0,1)$, let $B^{(t+1)}=B'$, if
    \[U \le \min\Big\{1, \frac{\pi^n(B')q(B^{(t)} \mid B')}{\pi^n(B^{(t)})q(B'\mid B^{(t)})}\Big\};\]
    otherwise, let $B^{(t+1)}=B^{(t)}$.
    \item Given $B^{(t+1)}$, obtain $\beta^{B^{(t+1)}}=(\beta^{B^{(t+1)}}_1,\ldots,\beta^{B^{(t+1)}}_{|B|})^\top$ via sampling 
    \[ \beta^{B^{(t+1)}}_s \sim  \N_{K+1}(\hat{\beta}_s^{B^{(t+1)}}, \tfrac{\sigma^2 v}{\sigma^2 + \alpha v}\{ Z_{B^{(t+1)}(s)}^\top Z_{B^{(t+1)}(s)}\}^{-1}), \]
    and then set $\theta^{B^{(t+1)}}=Z^{B^{(t+1)}}\beta^{B^{(t+1)}}$.
\end{enumerate}
Repeating this process $M$ times and discarding the first $m$ burn-in iterations, we obtain a sample of $(B^{(m+1)}, \theta^{B^{(m+1)}} )$, \ldots , $(B^{(M)}, \theta^{B^{(M)}} )$ from the joint posterior $\pi^n(B, \theta^B)$.  Then posterior mean of $\theta$ can be approximated by $(M-m)^{-1} \sum_{i=m+1}^M \theta^{B^{(i)}}$. Credible sets for any real-valued function $g(\theta)$ of $\theta$ can be obtained by obtaining quantiles of the samples 
\[ g(\theta^{B^{(i)}}), \quad i=m+1,\ldots,M. \]
For block configuration recovery, the maximum a posteriori (MAP) estimator for $B$ can be readily found by evaluating $\pi^n$, up to the normalizing constant, using formula \eqref{eq:marginal_B}, for each Monte Carlo sample $B^{(i)}$ and returning the maximizer.  For simplicity, we use a symmetric proposal distribution $q(B'\mid B)$ for the above algorithm, i.e., in each iteration, there is $0.5$ probability for a ``jump location'' to vanish and $0.5$ probability for a non-jump location to become a ``jump".

\section{Simulated data examples}
\label{S:examples}

\subsection{Methods}

In this section, from a perspective of numerical performance, we compare our proposed method to adaptive piecewise polynomial trend filtering in \citet{tibshirani2014adaptive}. We make use of the R package {\tt genlasso} for our implementation of trend filtering and the tuning parameter is chosen via five-fold cross-validation or the ``one-standard error'' rule, see \citet[Chapter~7]{hastie2009elements}. 

In order to implement the sampling procedures described above, some additional hyperparameters in \eqref{eq:marginal_B} also need to be specified. As mentioned before, since $\alpha=0.99$ has little practical difference from the $\alpha=1$ case, which corresponds to the genuine Bayesian model, we plug $\alpha=0.99$ into the posterior distribution functions for practical implementation. Next, for model variance $\sigma^2$, although the theory in Section~\ref{S:rate} assumes it to be known, in real applications, $\sigma^2$ may not be known and, therefore, must be estimated.  Of course, one can take a prior for $\sigma^2$ and get a corresponding joint posterior for $(\theta, \sigma^2)$; see \citet{martin2019empiricalpredict}.  Here, in keeping with the spirit of our empirical Bayes approach, we opt for a plug-in estimator.  Specifically, we consider
\begin{equation}
\label{eq:variance}
\hat\sigma^2 = \frac{1}{n} \sum_{i=1}^{n} (Y_{i} - \hat{\theta}^{\text{\sc tf}}_i)^2, 
\end{equation}
where $\hat{\theta}^{\text{\sc tf}}$ is the trend filtering/lasso estimate based on  cross-validation. For prior variance $v$,
it makes sense to take $v$ to be larger than $\sigma^2$ and, for the examples below,
with relatively small $\sigma^2$, we have found that $v = 1$ works well. Finally, $\lambda$ controls the penalty against large $|B|$ and, in the examples considered here, we conduct a sensitivity analysis in which $\lambda=0.2, 0.5, 1$ are all considered. 

For every data set, $5\times10^4$ iterations of the aforementioned MCMC algorithm, with an additional $5\times 10^4$ burn-in runs, are used to generate posterior samples.

\subsection{Scenarios}
\begin{figure}
\begin{center}
\begin{subfigure}[b]{0.5\textwidth}
\centering
\includegraphics[width=0.9\textwidth]{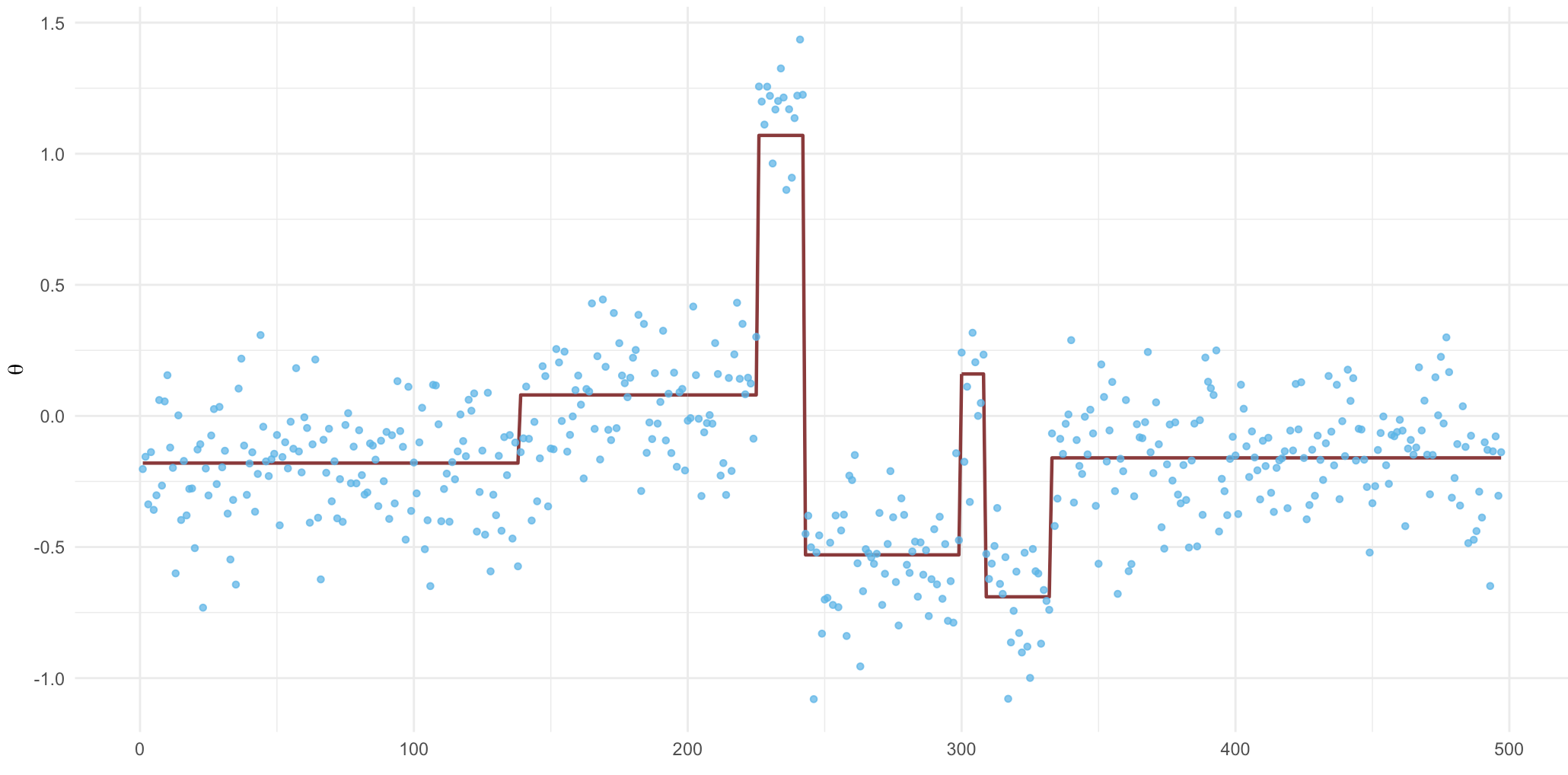}
\caption{Model~\ref{model:model1}}
\end{subfigure}%
\begin{subfigure}[b]{0.5\textwidth}
\centering
\includegraphics[width=0.9\textwidth]{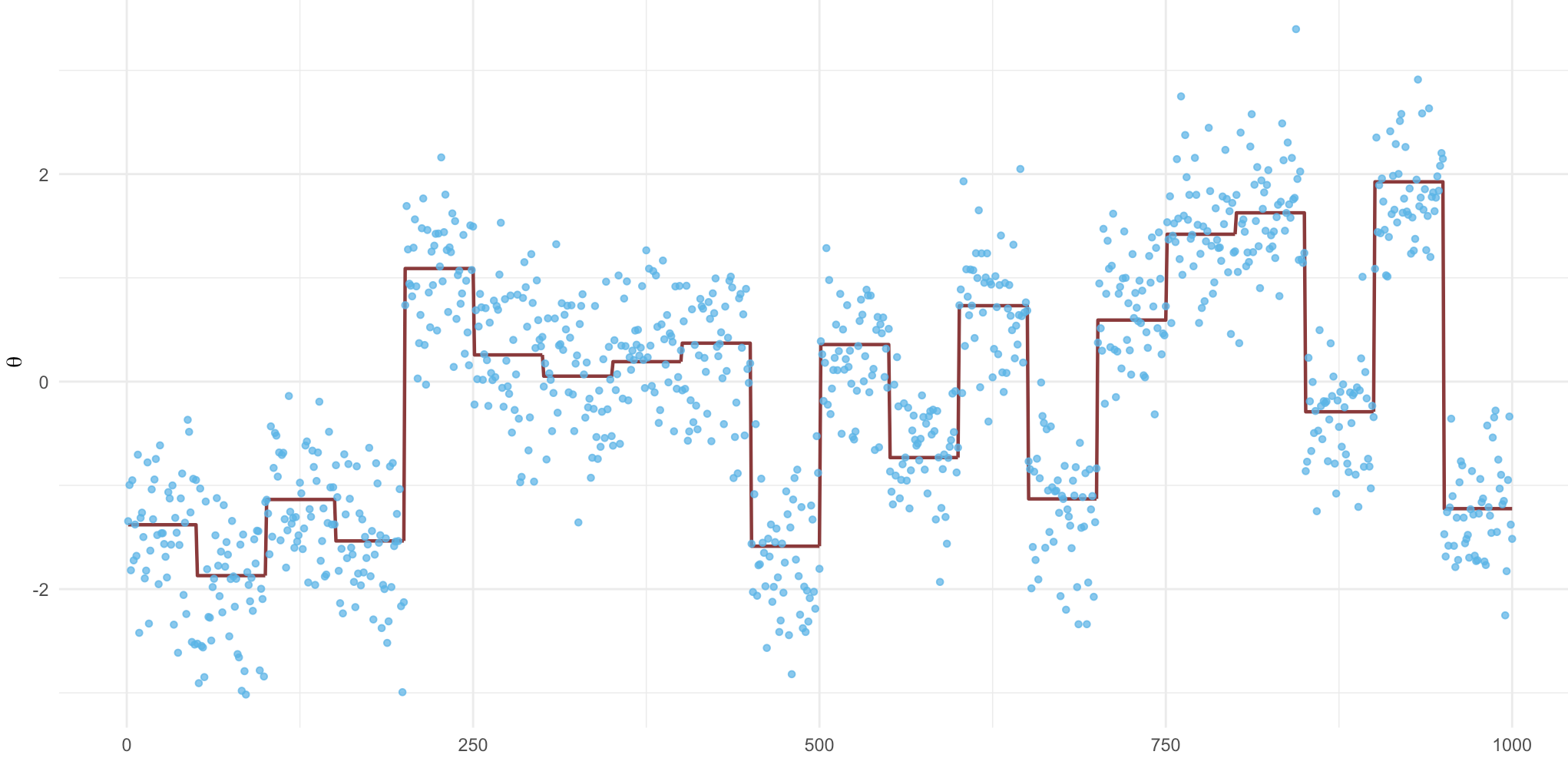}
\caption{Model~\ref{model:model2}}
\end{subfigure}%
\\ 
\begin{subfigure}[b]{0.5\textwidth}
\centering
\includegraphics[width=0.9\textwidth]{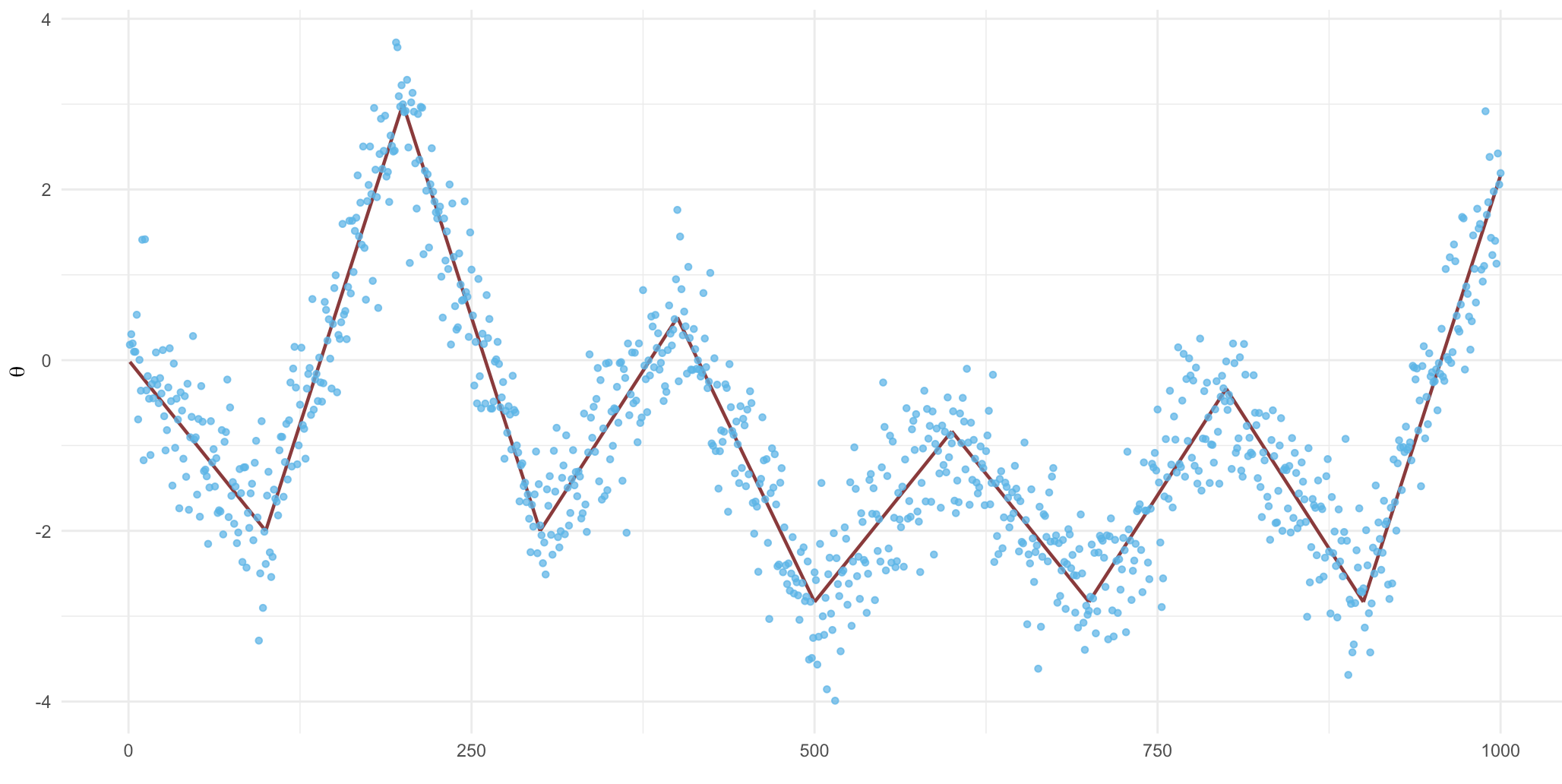}
\caption{Model~\ref{model:model3}}
\end{subfigure}%
\begin{subfigure}[b]{0.5\textwidth}
\centering
\includegraphics[width=0.9\textwidth]{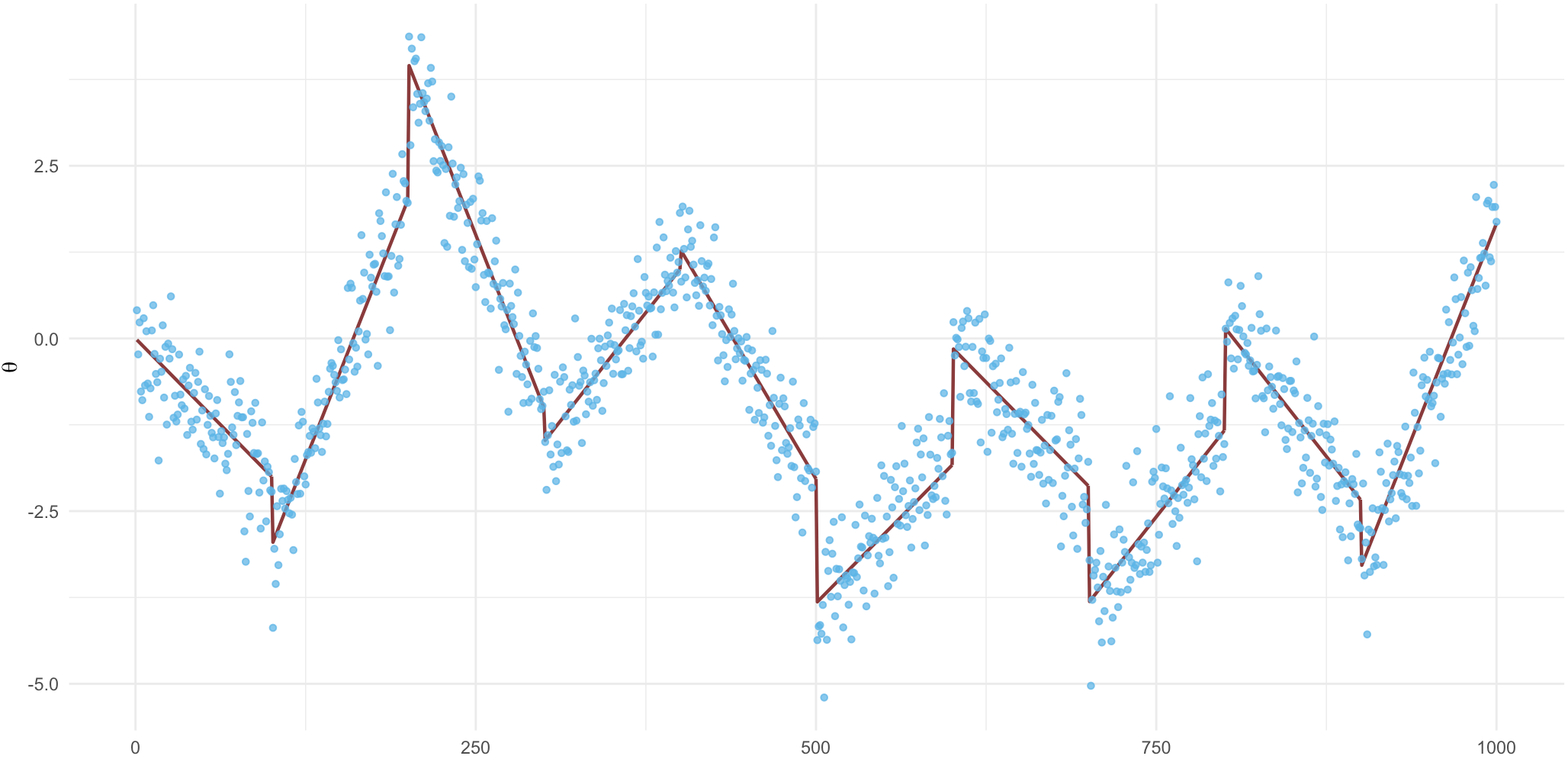}
\caption{Model~\ref{model:model4}}
\end{subfigure}%
\\
\begin{subfigure}[b]{0.5\textwidth}
\centering
\includegraphics[width=0.9\textwidth]{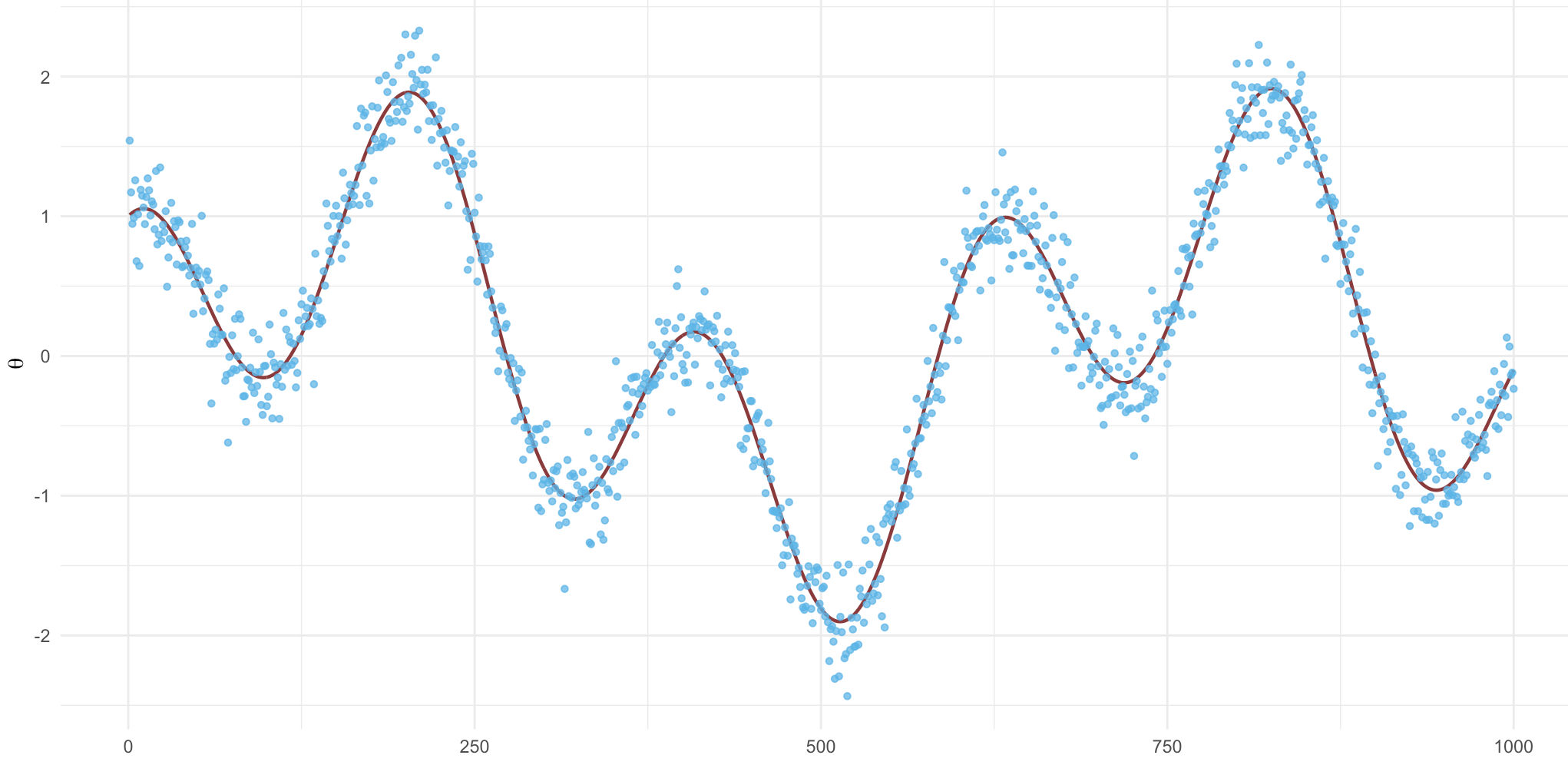}
\caption{Model~\ref{model:model5}}
\end{subfigure}%
\begin{subfigure}[b]{0.5\textwidth}
\centering
\includegraphics[width=0.9\textwidth]{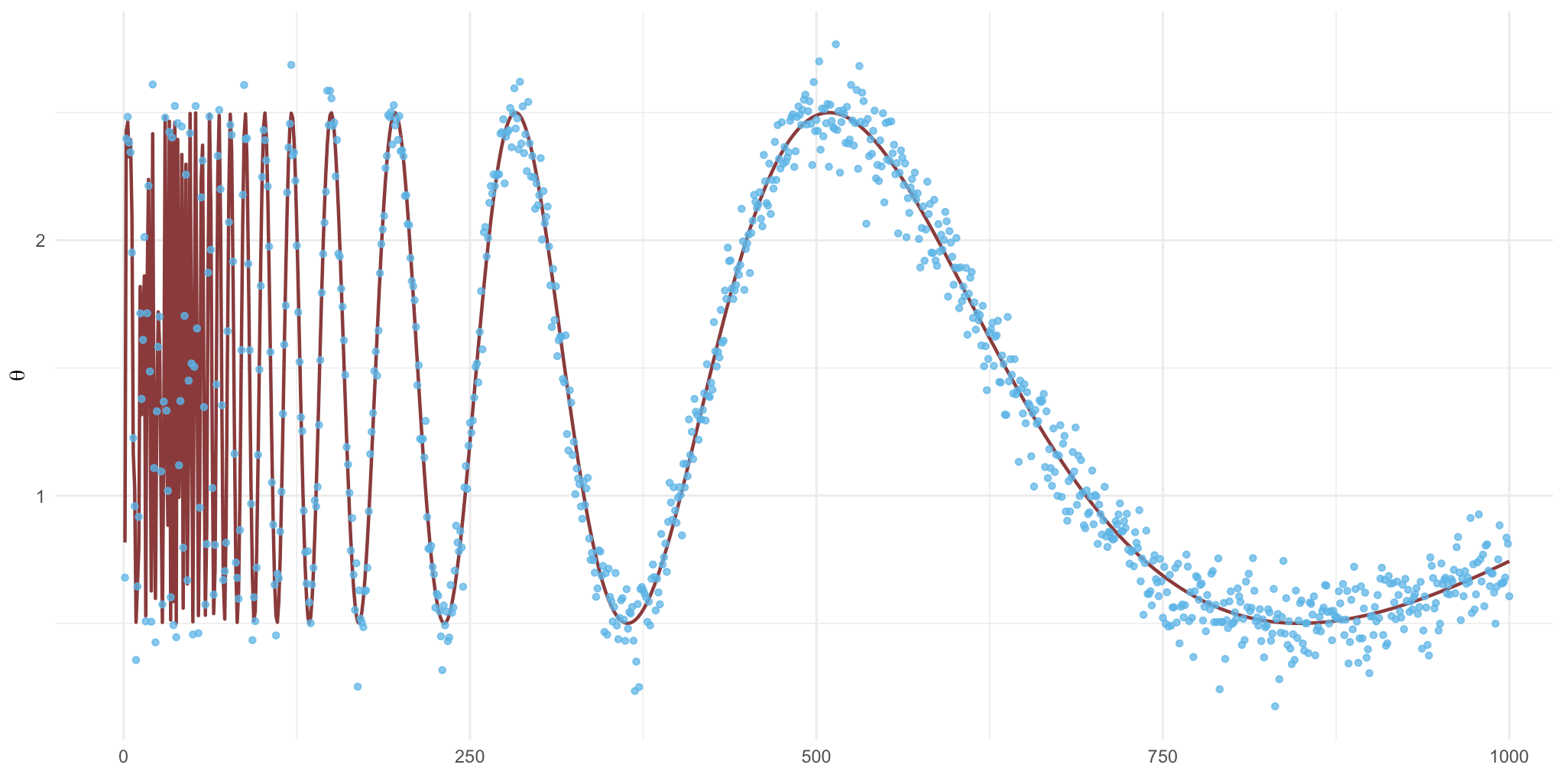}
\caption{Model~\ref{model:model6}}
\end{subfigure}%
\end{center}
\caption{Plots of the signal and a representative data set from each of the six models.}
\label{fig:models}
\end{figure}

For generating data, we consider the following six different models for the true signal $\theta^\star$.  The underlying truth and the simulated data are depicted in Figure~\ref{fig:models}.  

\begin{model}
\label{model:model1}
Piecewise constant model:  
\[\theta_i^\star=f_1(i), \quad \text{for} \quad i=1,\ldots,n, \quad n=497,\]
where $f_1(\cdot)$ is a piecewise constant function with $|B^*|=7$ as in \citet{frick2014multiscale}, p.~561; see, also, \citet{fryzlewicz2014wild}, Appendix~B(2). See Figure~\ref{fig:models}(a). The data is generated by
\[Y_i \sim \N(\theta_i^\star, 0.04), \quad i=1,\ldots,n, \quad n=497.\]
\end{model}


\begin{model}
\label{model:model2}
Piecewise constant model:  
\[\theta_i^\star=f_2(i), \quad \text{for} \quad i=1,\ldots,n, \quad n=1000,\]
where $f_2(\cdot)$ is a piecewise constant function with $|B^*|=20$ as in \cite{fan2018approximate}; see Figure~\ref{fig:models}(b). The the data is generated by
\[Y_i \sim \N(\theta_i^\star, 0.25), \quad i=1,\ldots,n, \quad n=1000.\]
\end{model}


\begin{model}
\label{model:model3}
Piecewise linear model with continuous mean:  
\[\theta_i^\star=f_3(i), \quad \text{for} \quad i=1,\ldots,n, \quad n=1000,\]
where $f_3(\cdot)$ is a piecewise linear function with continuous means (signal only with a change in its slope); $|B^*|=10$; see Figure~\ref{fig:models}(c). The the data is generated by
\[Y_i \sim \N(\theta_i^\star, 0.25), \quad i=1,\ldots,n, \quad n=1000.\] 
\end{model}


\begin{model}
\label{model:model4}
Piecewise linear model with jumps in mean:  
\[\theta_i^\star=f_4(i), \quad \text{for} \quad i=1,\ldots,n, \quad n=1000,\]
where $f_4(\cdot)$ is a piecewise linear function with signal both changing in its slope and intercept; $|B^*|=10$; see Figure~\ref{fig:models}(d). The the data is generated by
\[Y_i \sim \N(\theta_i^\star, 0.25), \quad i=1,\ldots,n, \quad n=1000.\] 
\end{model}


\begin{model}
\label{model:model5}
Trigonometric wave:  
\[\theta_i^\star=\sin(i/100)+\cos(i/33), \quad \text{for} \quad i=1,\ldots,n, \quad n=1000,\]
see Figure~\ref{fig:models}(e). The the data is generated by
\[Y_i \sim \N(\theta_i^\star, 0.04), \quad i=1,\ldots,n, \quad n=1000.\] 
\end{model}


\begin{model}
\label{model:model6}
Doppler wave as in \citet{tibshirani2014adaptive}, p.~298:
\[\theta_i^\star=\sin(4n/i)+1.5, \quad \text{for} \quad i=1,\ldots,n, \quad n=1000,\]
 see Figure~\ref{fig:models}(f). The the data is generated by
\[Y_i \sim \N(\theta_i^\star, 0.04), \quad i=1,\ldots,n, \quad n=1000.\] 
\end{model}


\subsection{Results}
\begin{figure}[t]
\begin{center}
\begin{subfigure}[b]{0.5\textwidth}
\centering
\includegraphics[width=0.99\textwidth]{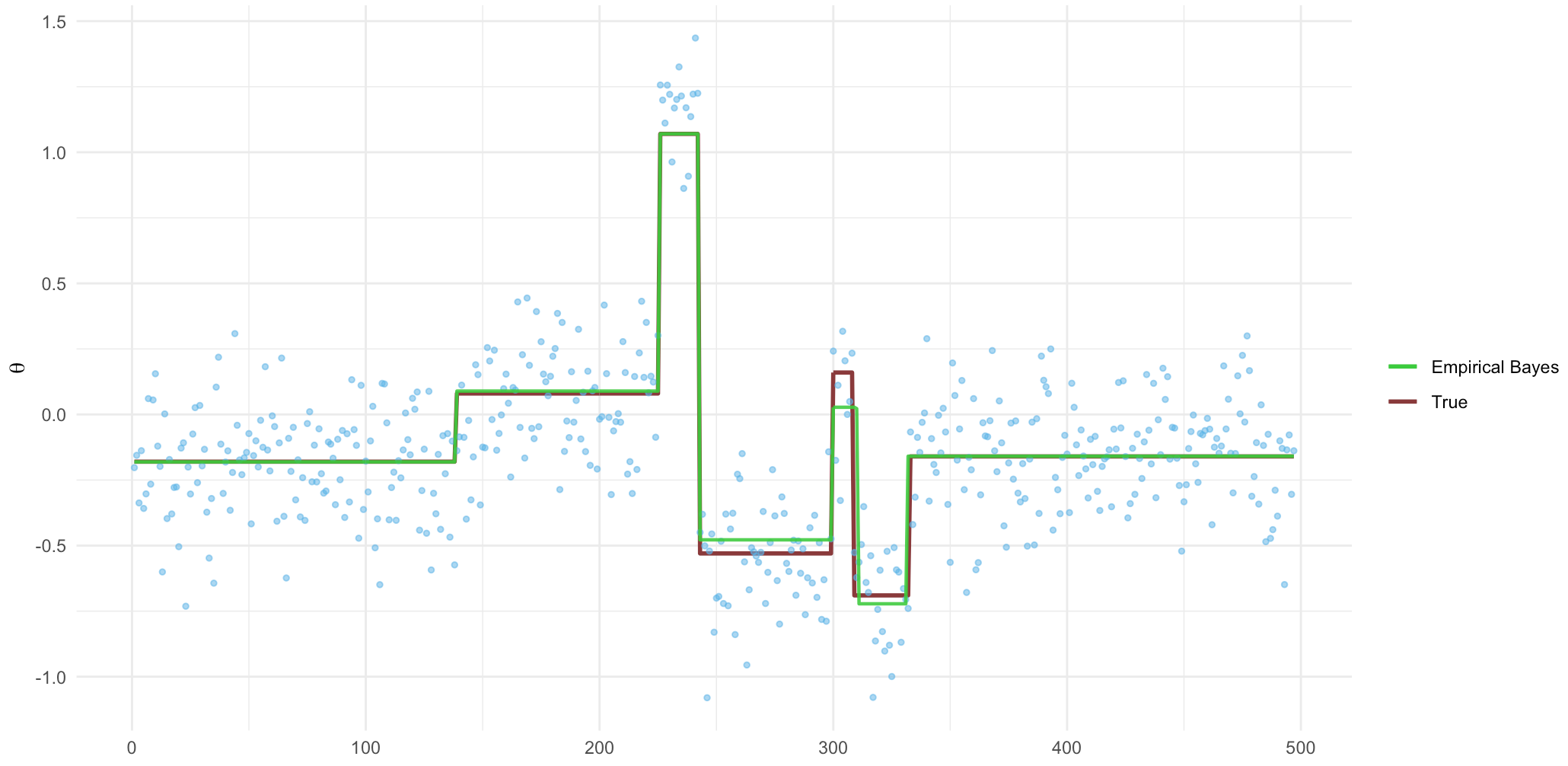}
\caption{Model~\ref{model:model1}: Empirical Bayes}
\end{subfigure}%
\begin{subfigure}[b]{0.5\textwidth}
\centering
\includegraphics[width=0.99\textwidth]{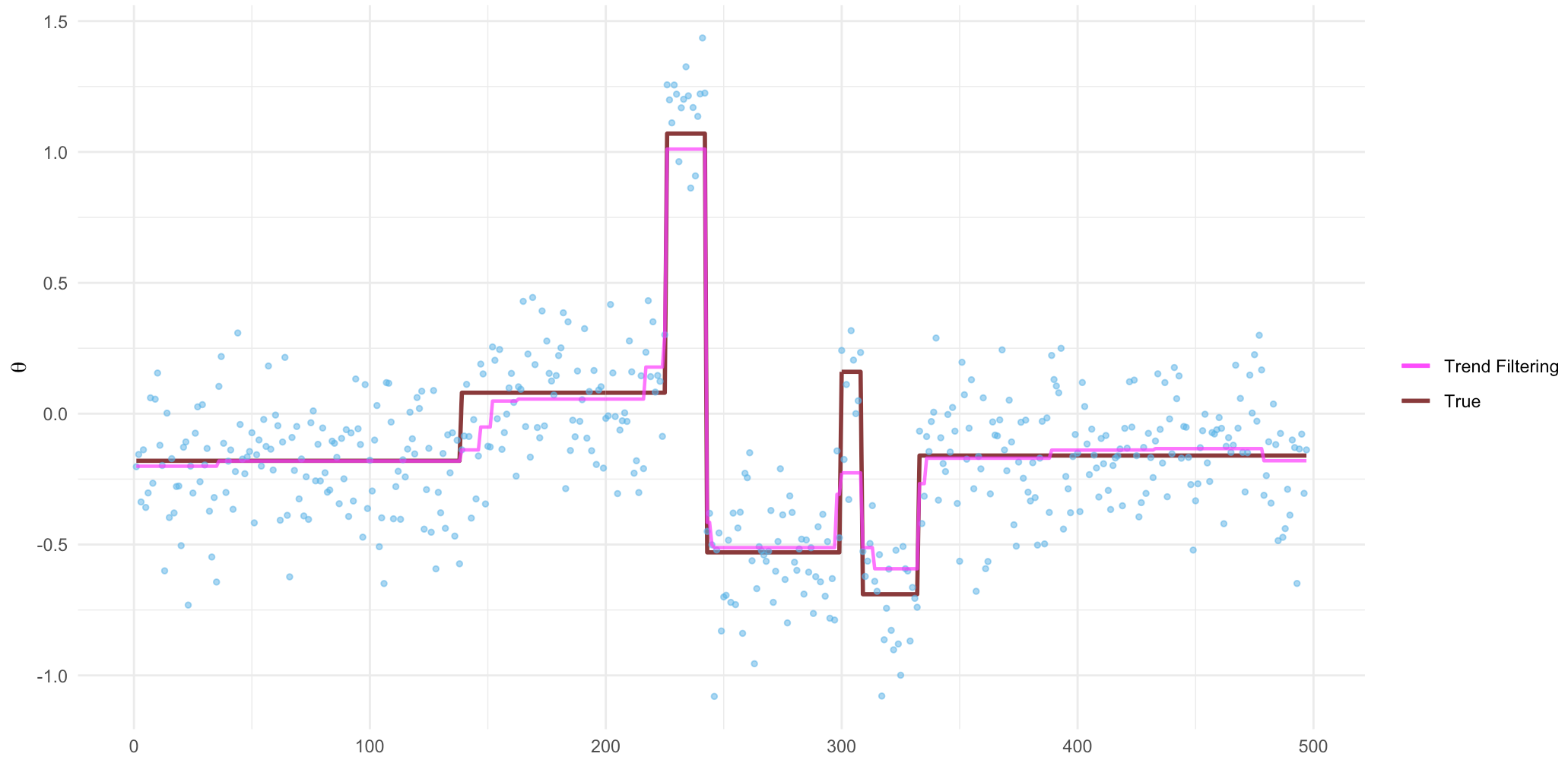}
\caption{Model~\ref{model:model1}: Trend filtering}
\end{subfigure}%
\\
\begin{subfigure}[b]{0.5\textwidth}
\centering
\includegraphics[width=0.99\textwidth]{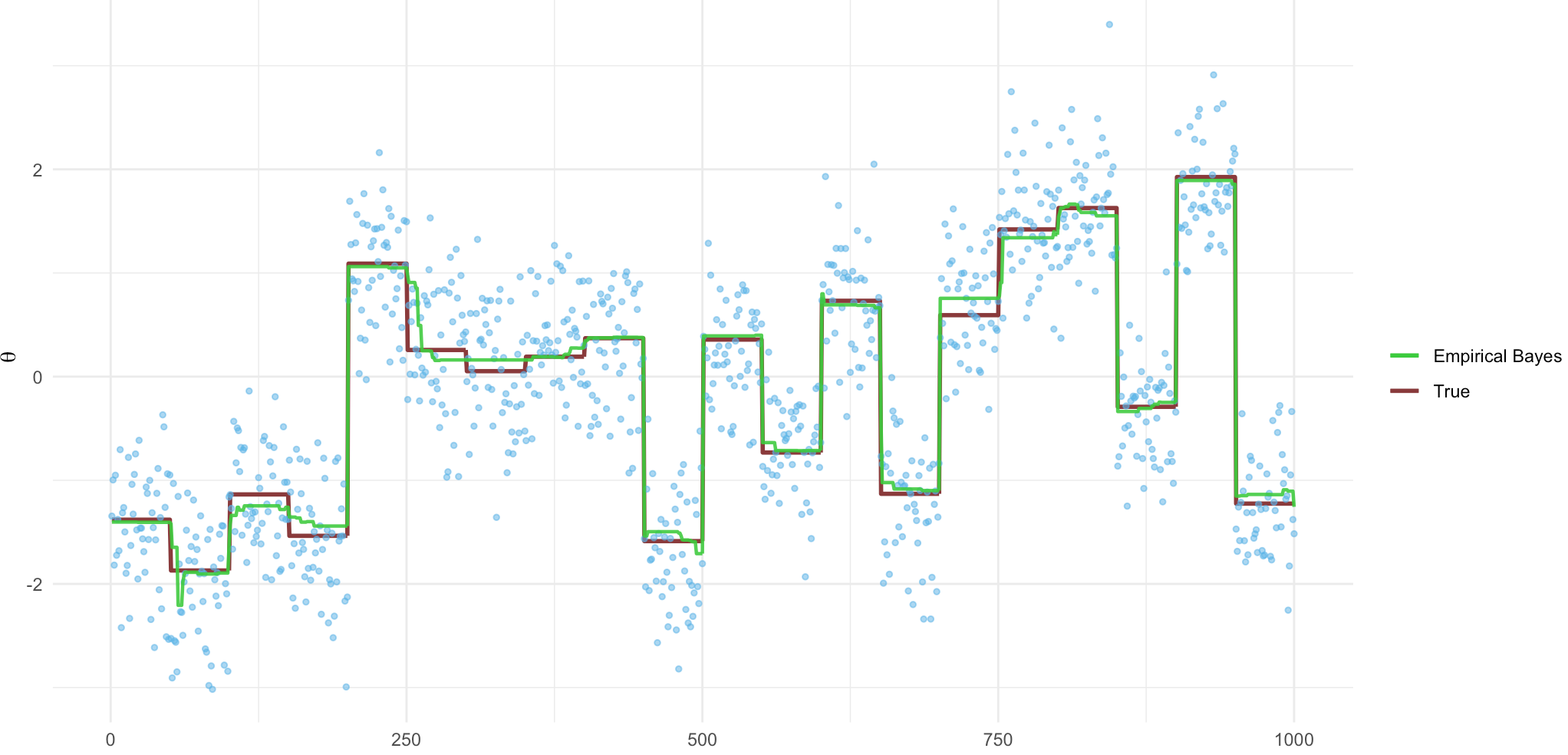}
\caption{Model~\ref{model:model2}: Empirical Bayes}
\end{subfigure}%
\begin{subfigure}[b]{0.5\textwidth}
\centering
\includegraphics[width=0.99\textwidth]{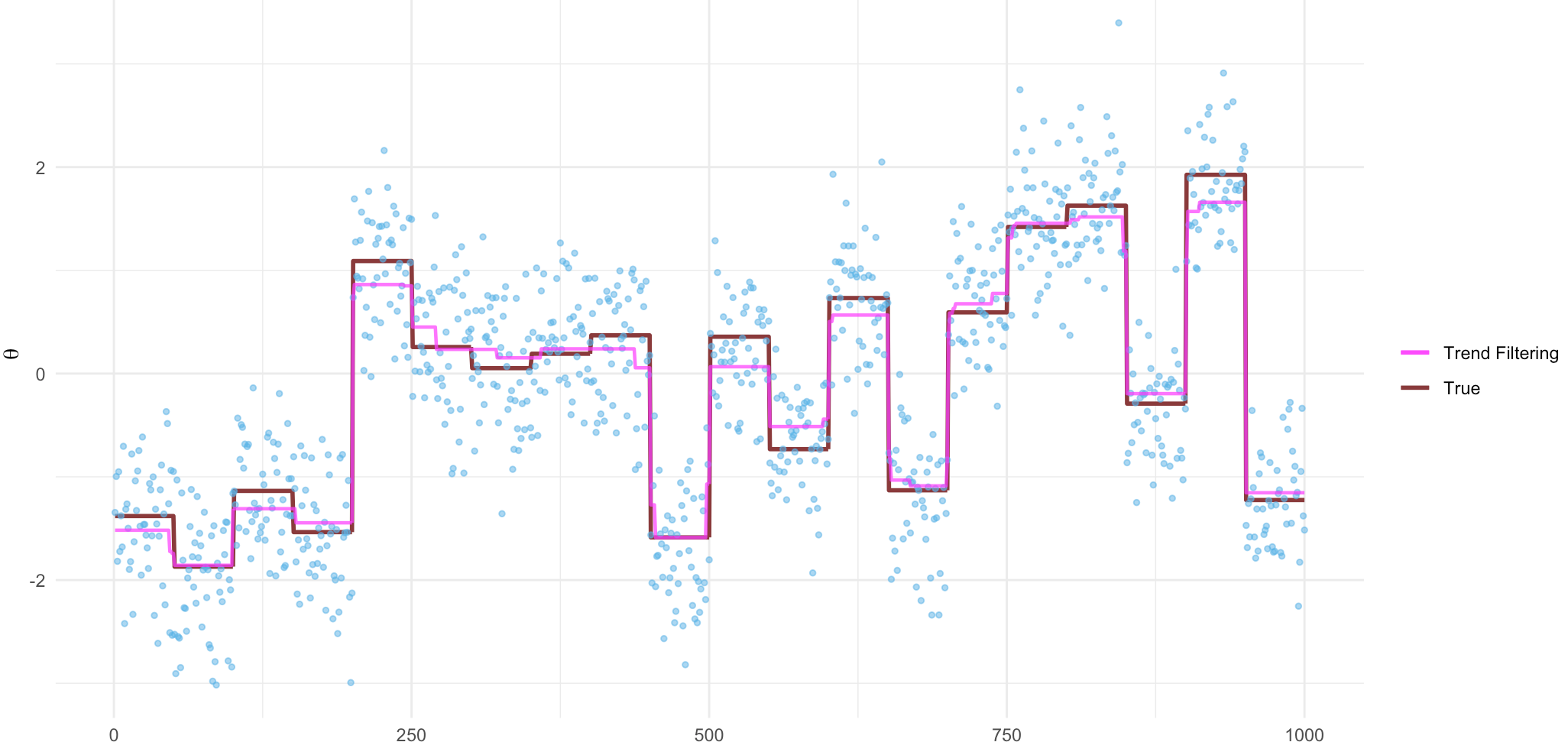}
\caption{Model~\ref{model:model2}: Trend filtering}
\end{subfigure}%
\\
\begin{subfigure}[b]{0.5\textwidth}
\centering
\includegraphics[width=0.99\textwidth]{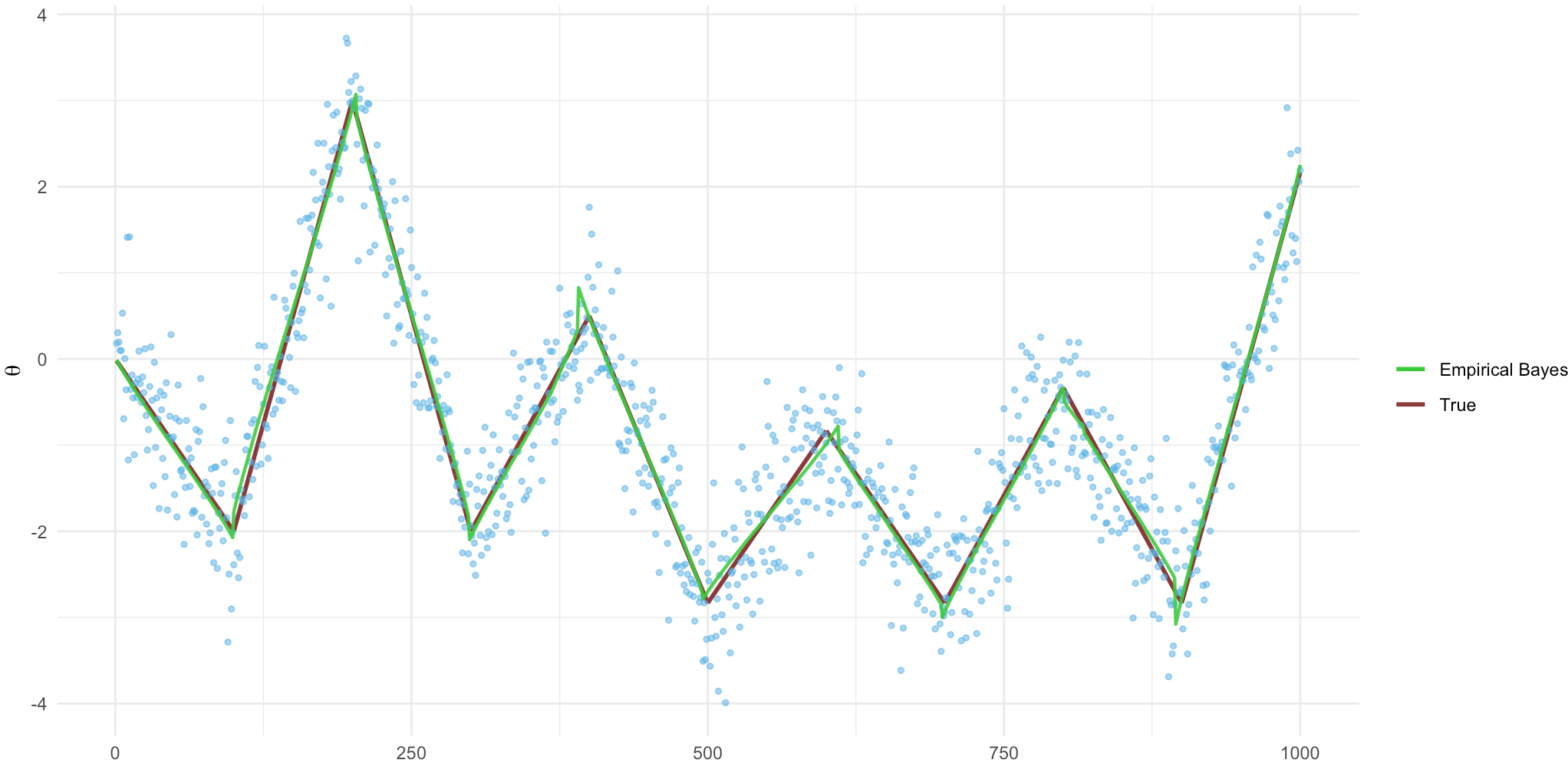}
\caption{Model~\ref{model:model3}: Empirical Bayes}
\end{subfigure}%
\begin{subfigure}[b]{0.5\textwidth}
\centering
\includegraphics[width=0.99\textwidth]{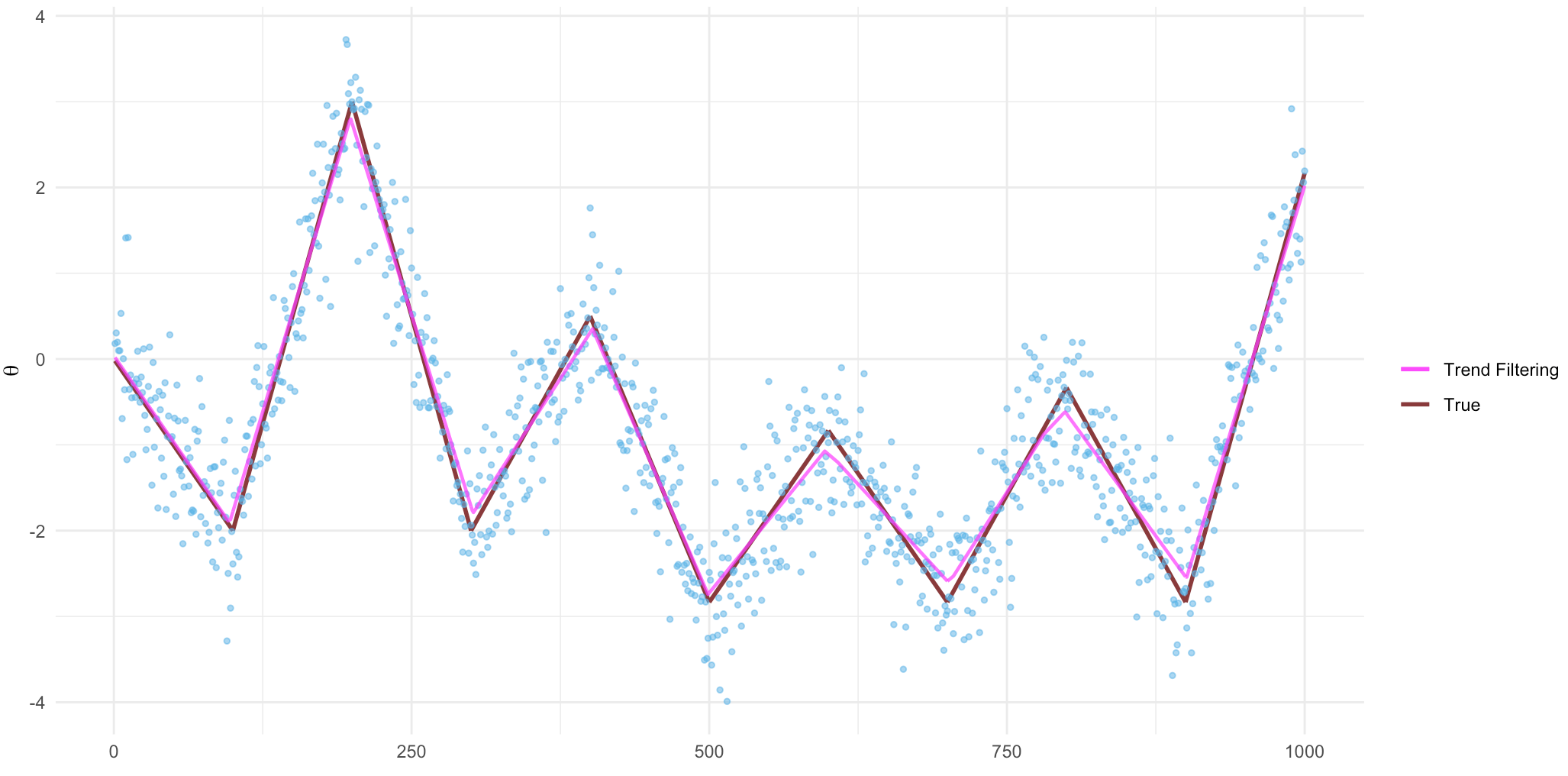}
\caption{Model~\ref{model:model3}: Trend filtering}
\end{subfigure}%
\end{center}
\caption{Plots of the empirical Bayes and trend filtering estimates of the signal for representative cased under Models~\ref{model:model1}--\ref{model:model3}.}
\label{fig:estimates.a}
\end{figure}

\begin{figure}[t]
\begin{center}
\begin{subfigure}[b]{0.5\textwidth}
\centering
\includegraphics[width=0.99\textwidth]{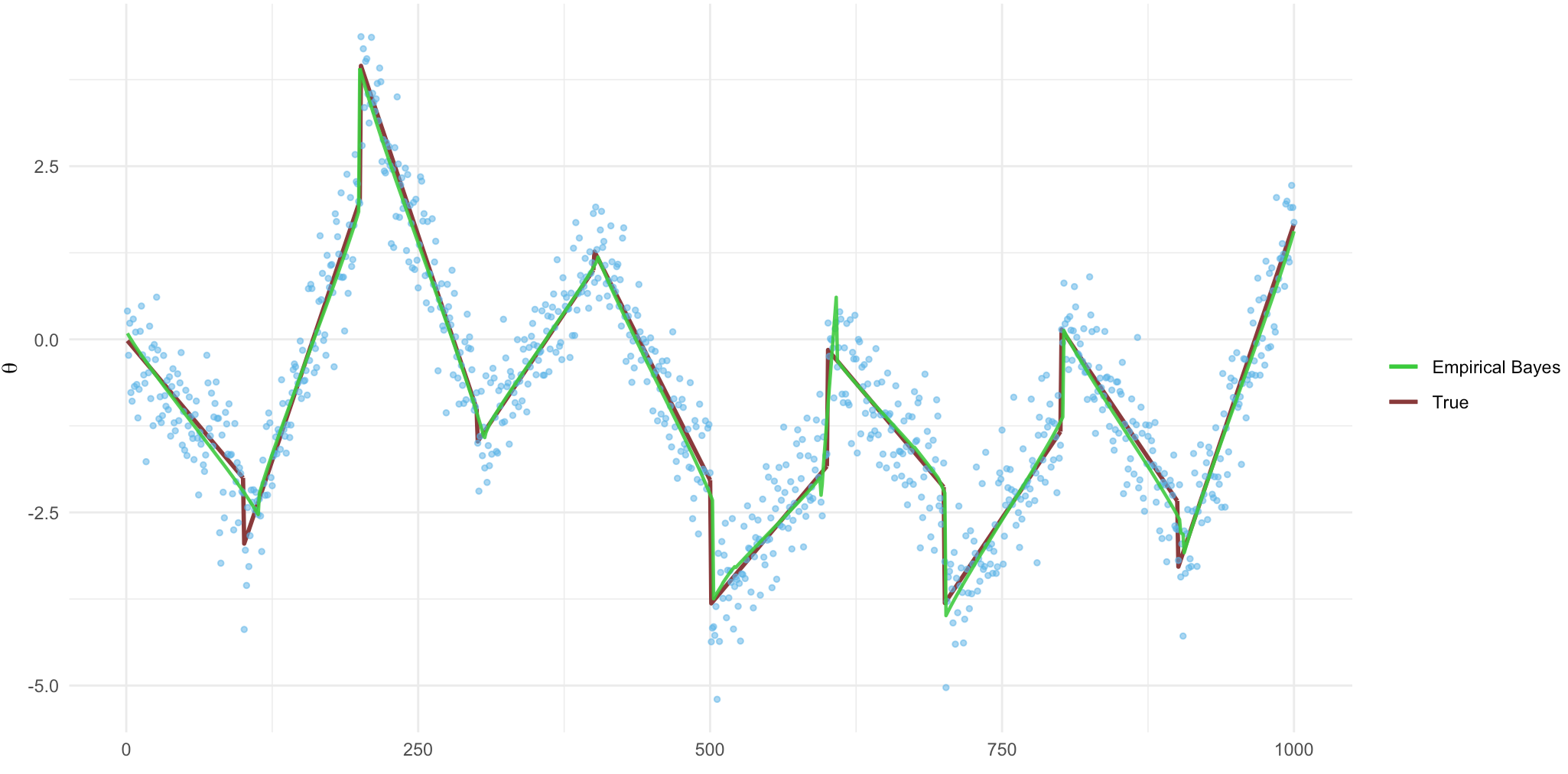}
\caption{Model~\ref{model:model4}: Empirical Bayes}
\end{subfigure}%
\begin{subfigure}[b]{0.5\textwidth}
\centering
\includegraphics[width=0.99\textwidth]{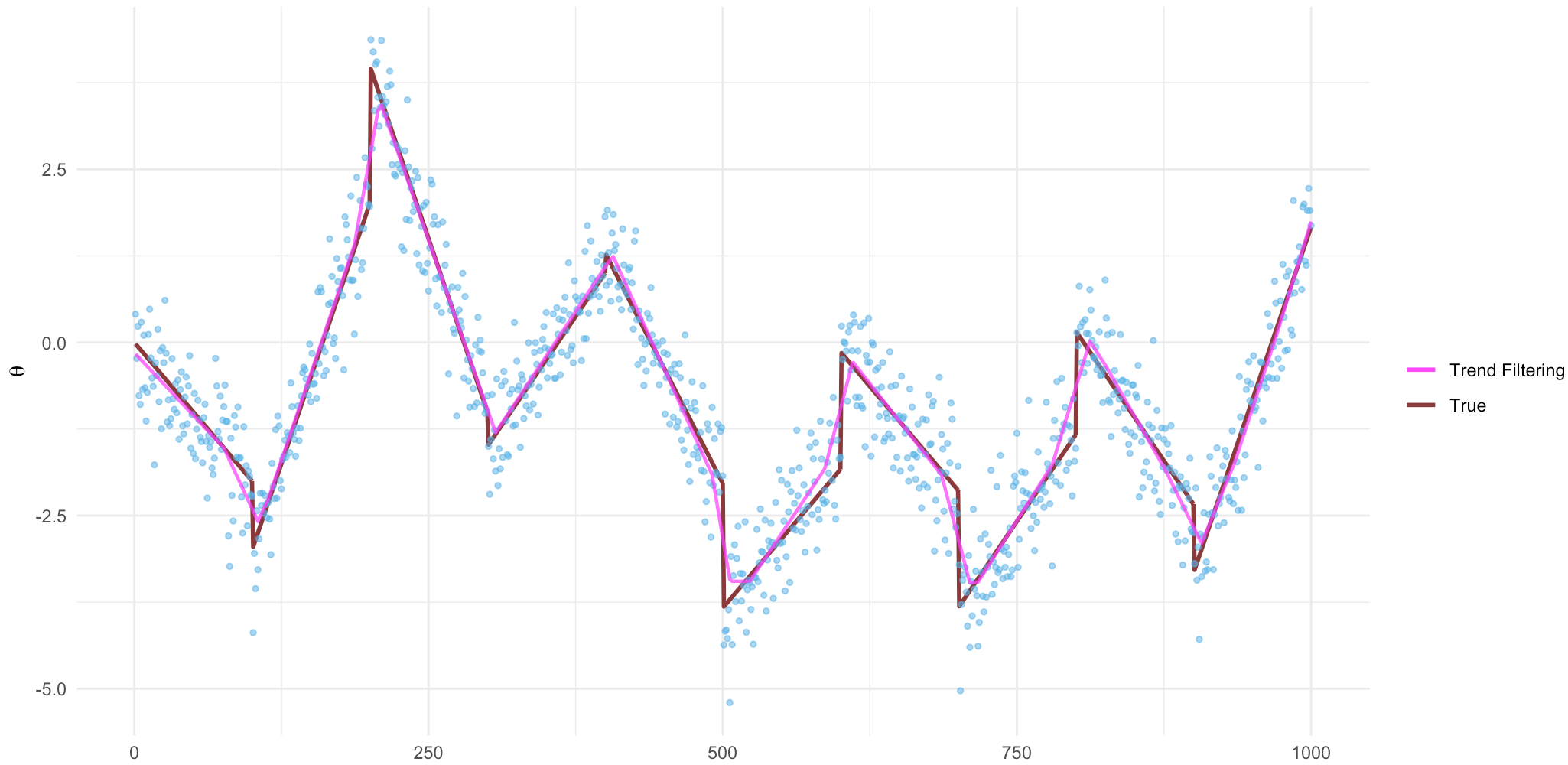}
\caption{Model~\ref{model:model4}: Trend filtering}
\end{subfigure}%
\\
\begin{subfigure}[b]{0.5\textwidth}
\centering
\includegraphics[width=0.99\textwidth]{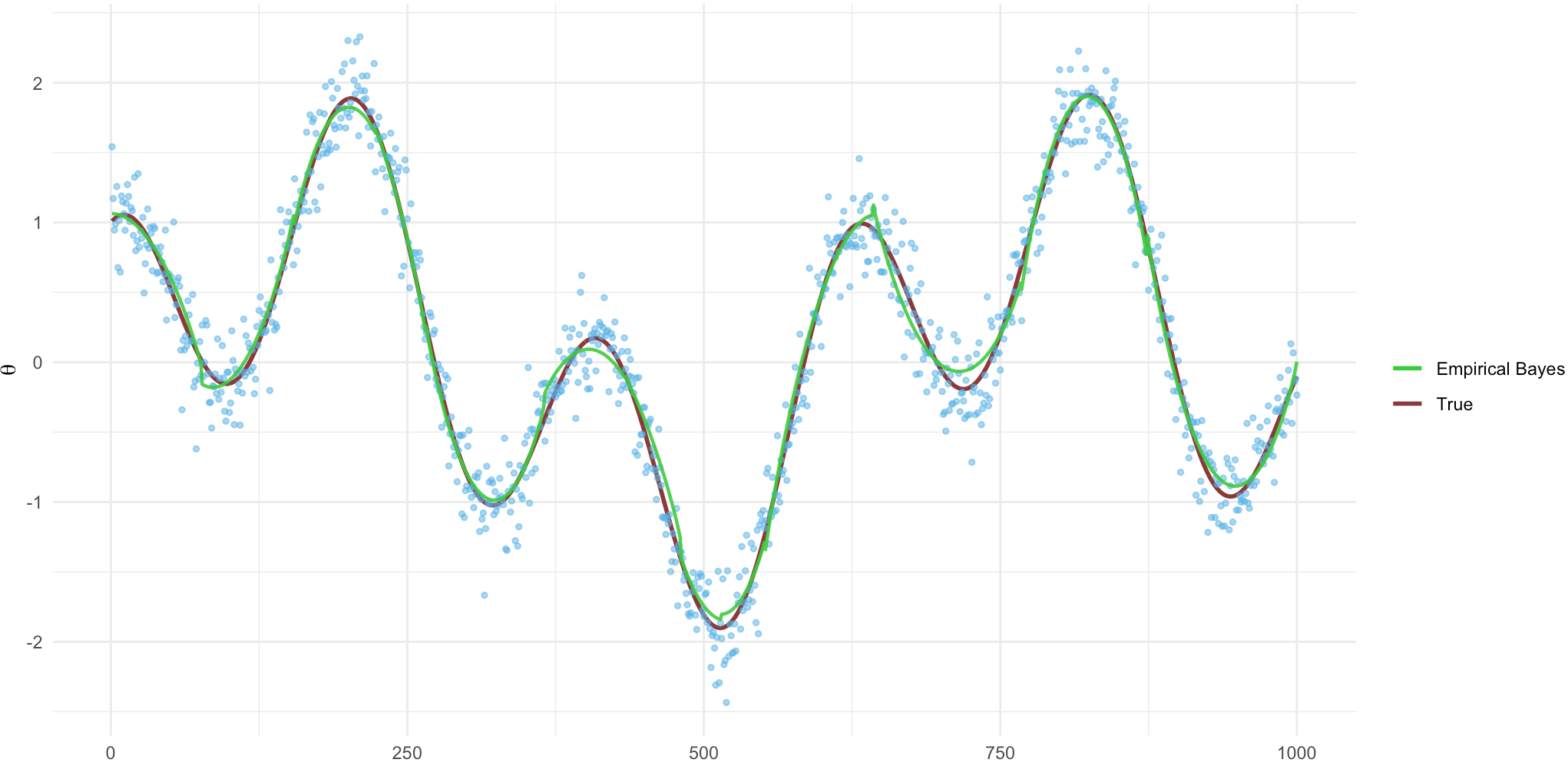}
\caption{Model~\ref{model:model5}: Empirical Bayes}
\end{subfigure}%
\begin{subfigure}[b]{0.5\textwidth}
\centering
\includegraphics[width=0.99\textwidth]{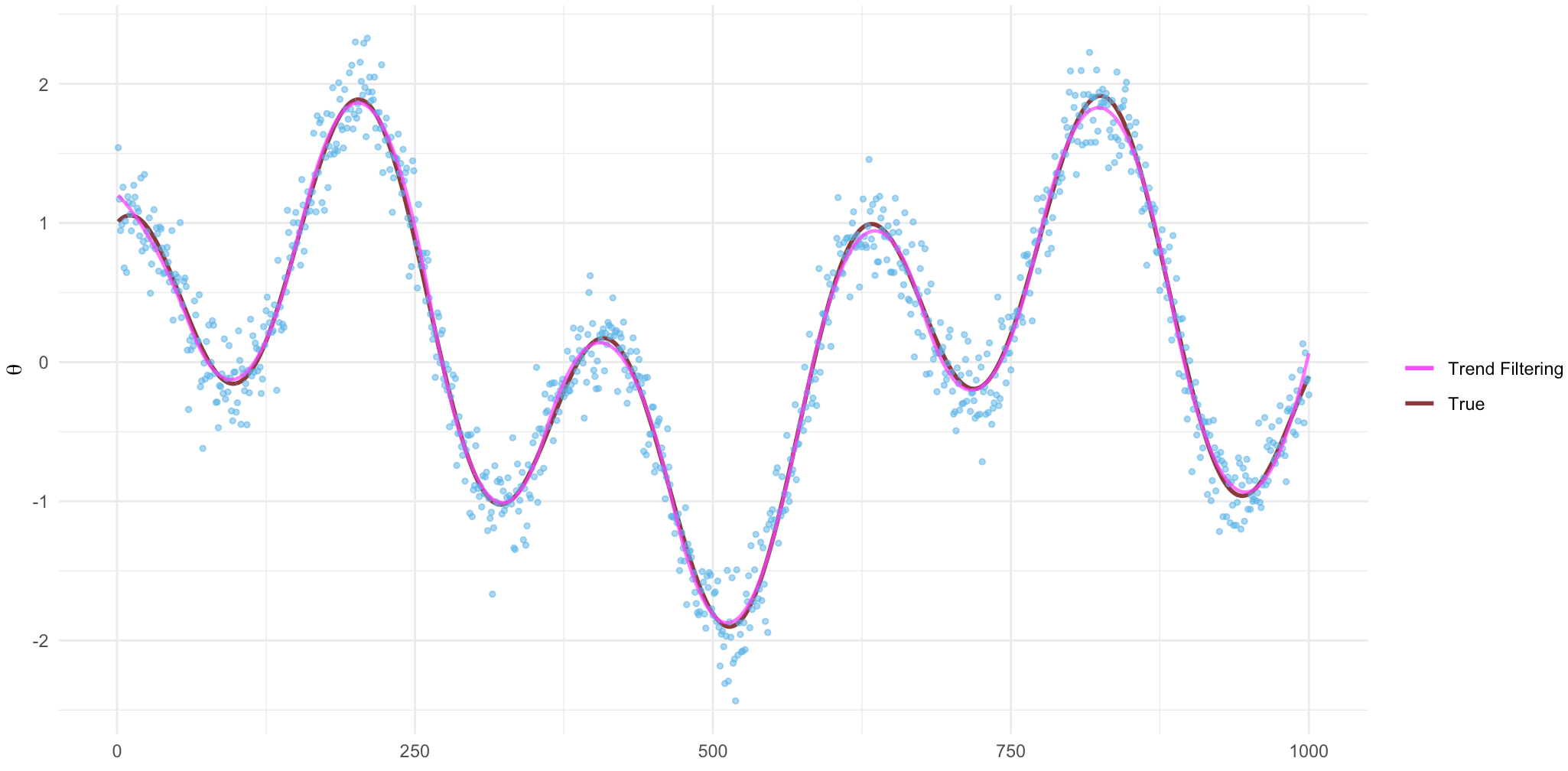}
\caption{Model~\ref{model:model5}: Trend filtering}
\end{subfigure}%
\\
\begin{subfigure}[b]{0.5\textwidth}
\centering
\includegraphics[width=0.99\textwidth]{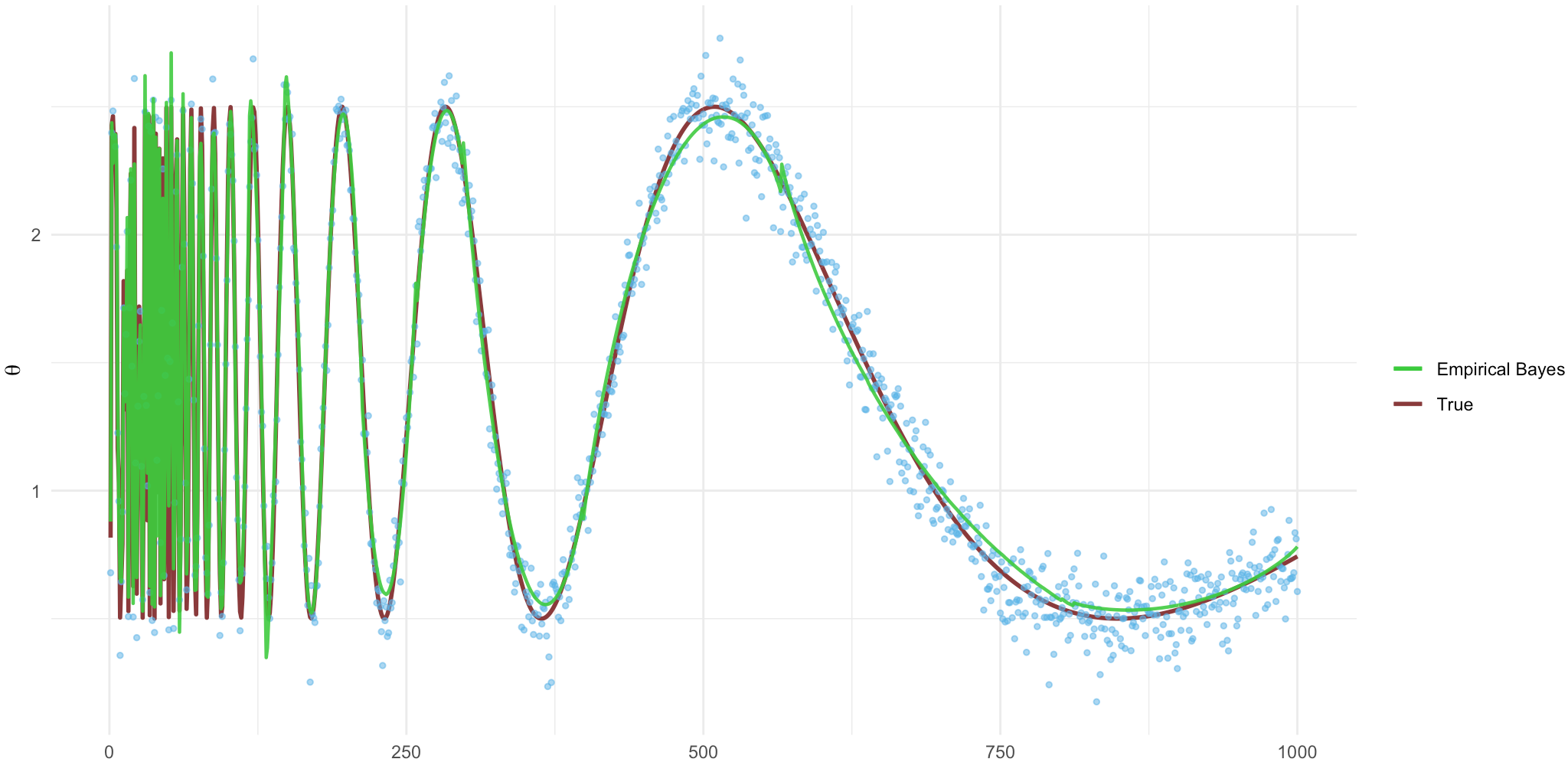}
\caption{Model~\ref{model:model6}: Empirical Bayes}
\end{subfigure}%
\begin{subfigure}[b]{0.5\textwidth}
\centering
\includegraphics[width=0.99\textwidth]{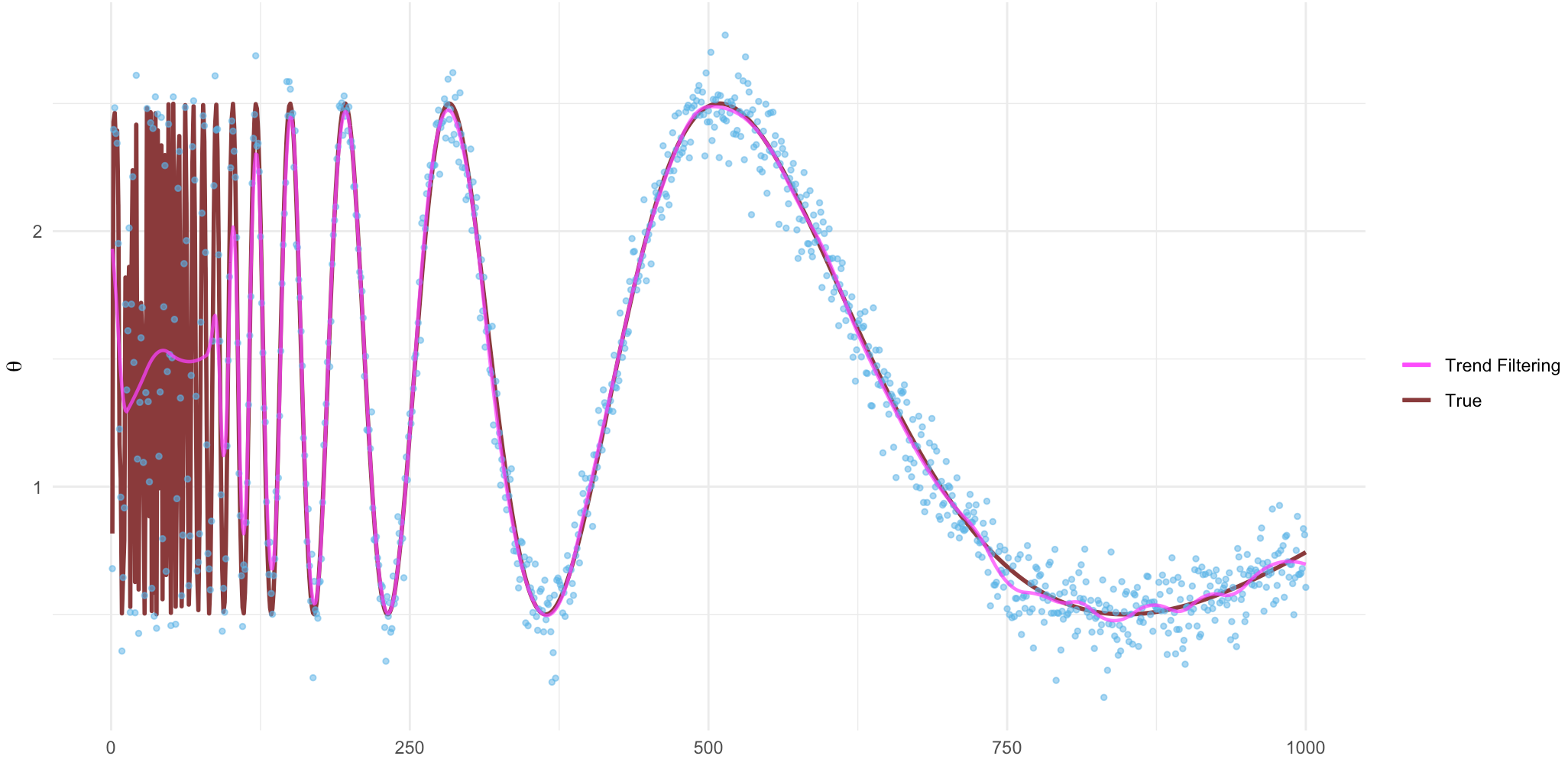}
\caption{Model~\ref{model:model6}: Trend filtering}
\end{subfigure}%
\end{center}
\caption{Plots of the empirical Bayes and trend filtering estimates of the signal for representative cased under Models~\ref{model:model4}--\ref{model:model6}.}
\label{fig:estimates.b}
\end{figure}

We investigate numerical performance of the two methods in terms of estimation error and block selection accuracy. For Models~\ref{model:model1}--\ref{model:model6}, squared estimation error loss, $\|\hat\theta-\theta^\star\|^2$, is computed in Table~\ref{tab:fit} where $\hat\theta$ is either our posterior mean or the trend filtering estimator obtained by cross-validation.  In addition, the estimated signal $\hat{\theta}$ and the true $\theta^\star$ are plotted in Figures~\ref{fig:estimates.a}--\ref{fig:estimates.b}. Note that for those graphical comparisons, the trend filtering estimator is computed through ``one-standard error'' rule, since it is usually smoother than that chosen by  cross-validation, although it typically suffers from higher mean squared error; see \citet[Chapter~7]{hastie2009elements} for details.

For structure recovery results, since it is more meaningful to discuss about change-points/block partitions in lower order piecewise polynomial models, such as piecewise constant and piecewise linear, here we only focus on  Models~\ref{model:model1}--\ref{model:model2} (results displayed in Tables~\ref{tab:model1}--\ref{tab:model2}).
Estimated block partition $\hat{B}$ of trend filtering is obtained by looking at the nonzero entries of $D^{(K+1)}\hat{\theta}$, i.e., the $K^{\text{th}}$-order ``knots" of $\hat{\theta}$; see \citet{guntuboyina2020adaptive}. For our empirical Bayes method, $\hat{B}$ is the maximizer of marginal posterior probability $\pi^n(B)$.
 To get a better understanding of structure learning performances, we consider using multiple criteria to measure the change-point detection/block selection accuracy.  Among the 100 replications for each model, we calculate the probability of  $\hat{B}$ equaling the true $B^\star$, covering the true $B^\star$, and being nested in the true $B^\star$, which are denoted as $P(\hat{B}=B^\star)$, $P(\hat{B}\supset B^\star)$ and $P(\hat{B}\subset B^\star)$ respectively. $|{\hat{B}}|$, the size of $\hat{B}$, is also reported. In addition, as is discussed in Section~\ref{SS:structure}, an equivalent representation of block partition is the set of jump locations $J$, defined in Theorem~\ref{thm:B_consist}. Hausdorff distance between $J$ and $J^\star$ can be calculated using the following formula,
 \[ H(J \mid J^\star) = \max_{j^\star \in J^\star} \min_{j \in J} |j - j^\star|+
 \max_{j \in J} \min_{j^\star \in J^\star} |j - j^\star|. \]
 Finally, we consider an $(n-K-1)$-dimensional binary vector $S$, with $S_i=1$ if and only if $i \in J$. Hamming distance between $\hat{S}$ and $S^\star$ is reported as a measure of how close the $\hat{J}$ and $J^\star$ are from each other.
 
 From the perspective of estimation accuracy for $\theta^\star$, in Table~\ref{tab:fit}, our method achieves smaller squared error loss compared to trend filtering except for Models~\ref{model:model3} and \ref{model:model5} which are the two that are continuous; the Doppler wave function in Model~\ref{model:model6} is continuous too, but the high frequency oscillation part in $[0, 100]$ makes it ``almost discontinuous.'' Therefore, our method tends to have an advantageous estimation performance when the underlying $\theta^\star$ have jump discontinuities, particularly for piecewise constant signals, which is discontinuous in its nature. Furthermore, our method demonstrates stronger structure recovery for piecewise constant Models~\ref{model:model1} and \ref{model:model2}. Compared to trend filtering which tends to select more blocks, when $\lambda=0.5$ our method detects the exact block number for Model~\ref{model:model1}. In terms of Hamming distance and Hausdorff distance, our method also outperforms trend filtering for both models. However, the probabilities of identifying, covering and being nested in the true block partition are low for both methods. It is likely that identifying the exact change points is actually a challenging problem, given the high-dimensionality, i.e., there are hundreds or thousands of candidate points to be considered as change points/jump points.

\begin{center}
 \begin{table}[t]
 \centering
 {
 \begin{tabular} {c  c c c c c c} 
 \hline
Method & Model 1 & Model 2& Model 3 & Model 4& Model 5& Model 6 \\ 
 \hline
 Empirical Bayes&  1.2179&  26.4856& 14.7871 & 20.0689 & 3.6054&  10.6416\\
 $\lambda=1.0$&(0.0666) &(0.7748) & (0.4761)& (0.5481)& (0.0923) &(0.2542)\\
 \hline
 Empirical Bayes& 0.8189 &   16.3953 &10.7919 &  \textbf{15.2974}& 2.8549 &8.0418  \\
 $\lambda=0.50$&(0.0474)  & (0.5768)& (0.3288) &(0.4148) & (0.0880) &(0.1523)\\
 \hline
 Empirical Bayes& \textbf{0.7024} &  \textbf{13.9412} & 13.0632 & \textbf{15.4447} & 2.8367&  \textbf{7.3547}\\
 $\lambda=0.20$& (0.0334) &(0.4996) &(0.4353) & (0.5401)&(0.0719) &(0.1350)\\
 \hline
 Trend Filtering & 1.5208 &20.0122  & \textbf{6.8481} &  29.9949 & \textbf{1.2644}&  45.2222\\
 (cross-validation)&(0.0384) & (0.3709)& (0.1805)&(0.2882) & (0.0381)&(0.7007)\\
 \hline
 \end{tabular}
 }
\caption{Squared error loss between  $\hat{\theta}$ and $\theta^\star$ across 100 replications. Tuning parameter of trend filtering is selected by 5-fold cross-validation.}
\label{tab:fit}
\end{table}
\end{center}

\begin{table}[t]
\centering
{\small
\begin{tabular}{ c|c c c c c c} 
\hline
Method & $P(\hat{B}=B^\star)$ & $P(\hat{B}\supset B^\star)$ & $P(\hat{B}\subset B^\star)$ & Hamming & Hausdorff & $|\hat B|$ \\
    \hline
Empirical Bayes  &   0.03&  0.03& 0.03 & 4.14 &  3.86 & \textbf{7.00}\\
$\lambda=1.0$ &   (0.00)& (0.00)& (0.00) & (0.19) &(0.40) & (0.00)\\
\hline
Empirical Bayes  &   0.11&  0.11& 0.11 & 3.13 &  2.73 & \textbf{7.00}\\
$\lambda=0.50$ &   (0.00)& (0.00)& (0.00) & (0.18) &(0.31)  &(0.00)\\
\hline
Empirical Bayes  &   \textbf{0.19}&  0.19& \textbf{0.19} & \textbf{2.39} &  \textbf{2.57} &\textbf{7.01} \\
$\lambda=0.20$ &   (0.00)& (0.00)& (0.00) & (0.15) &(0.41) & (0.01) \\
\hline
Trend Filtering & 0.05 & \textbf{0.28} & 0.07 & 3.27 & 17.99  & 8.77 \\
(one std error) & (0.00) & (0.00) & (0.00) & (0.17) & (3.53) & (0.13) \\
\hline
\end{tabular}
}
\caption{Structure recovery results for Model~\ref{model:model1} ($|B^\star|=7$) across 100 replications.   The tuning parameter of trend filtering is chosen by one-standard-error rule.} 
\label{tab:model1}
\end{table}

\begin{table}[t]
\centering
{\small 
\begin{tabular}{ c|c c c c c c} 
\hline
Method & $P(\hat{B}=B^\star)$ & $P(\hat{B}\supset B^\star)$ & $P(\hat{B}\subset B^\star)$ & Hamming & Hausdorff & $|\hat B|$ \\
    \hline
Empirical Bayes  &   0.00&  0.00& 0.00 & 17.02 &  60.51 & 15.24 \\
$\lambda=1.0$ &   (0.00)& (0.00)& (0.00) & ( 0.74) &(4.38) &(0.29) \\
\hline
Empirical Bayes  &   0.00&  0.00& \textbf{0.01} & 14.27 &  55.35 & 15.96 \\
$\lambda=0.50$ &   (0.00)& (0.00)& (0.00) & (0.61) &(2.06)  & (0.23)\\
\hline
Empirical Bayes  &   0.00&  0.00& 0.00 & \textbf{12.06} &  \textbf{50.72} & \textbf{17.72}\\
$\lambda=0.20$ &   (0.00)& (0.00)& (0.00) & (0.61) &(1.50) & (0.24) \\
\hline
Trend Filtering  &   0.00&  \textbf{0.01}& 0.00 & 25.20 &  50.82  & 33.45\\
(one std error) &   (0.00)& (0.00)& (0.00) & (0.68) &(1.87) & (0.75) \\
\hline
\end{tabular}
}
\caption{Structure recovery results for Model~\ref{model:model2} ($|B^\star|=20$) across 100 replications.  The tuning parameter of trend filtering is chosen by one-standard-error rule.}
\label{tab:model2}
\end{table}

\section{Real data examples}
\label{S:real}

\subsection{DNA copy number analysis}
\begin{figure}[t]
\centering
    \begin{subfigure}[b]{0.5\textwidth}
                \centering
                \includegraphics[width=0.99\textwidth]{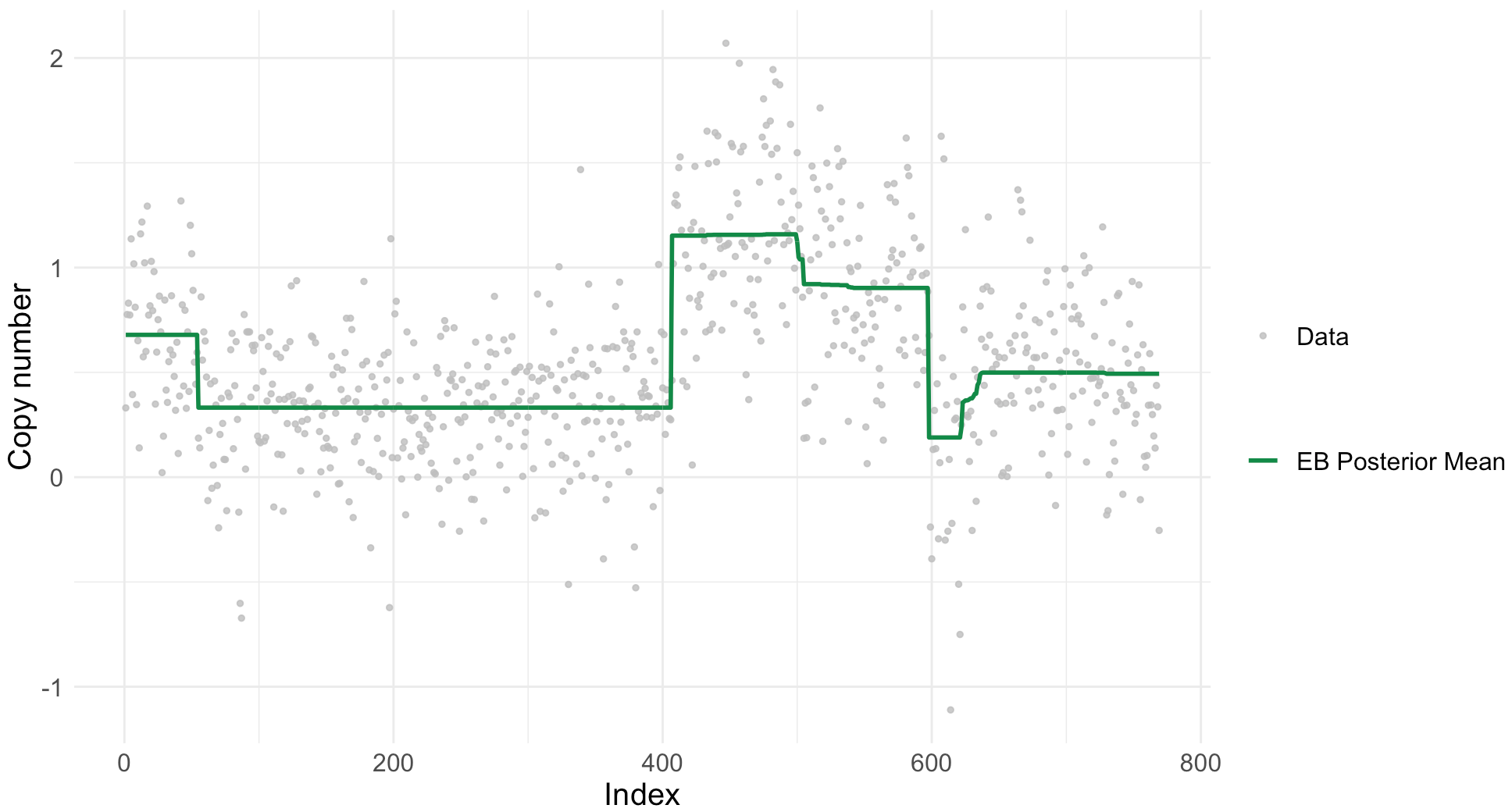}
            \caption{Posterior mean for $\theta$}
    \end{subfigure}%
    \begin{subfigure}[b]{0.5\textwidth}
                    \centering
                \includegraphics[width=0.99\textwidth]{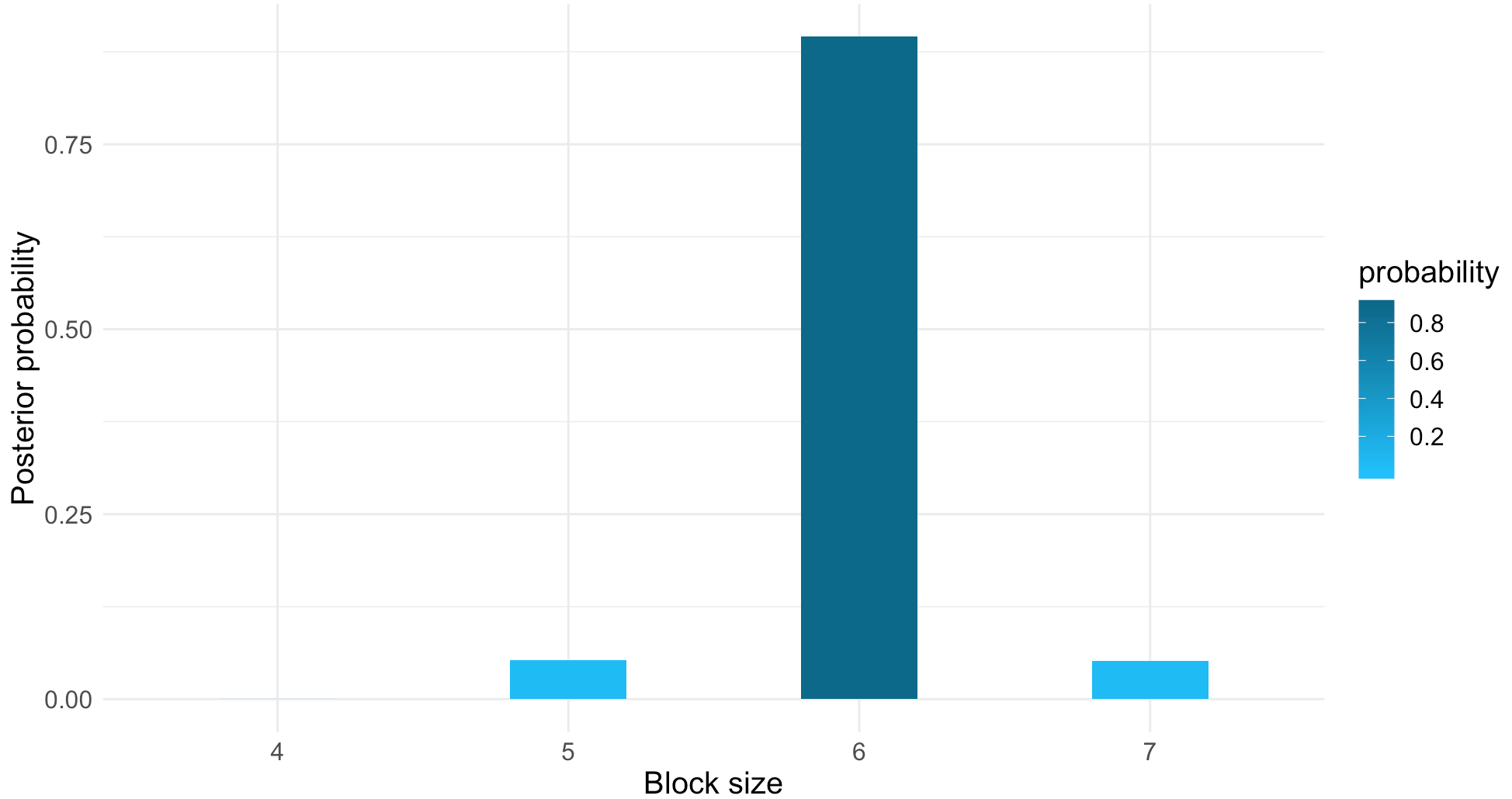}
         \caption{Posterior distribution for $|B|$}
    \end{subfigure}%
\caption{DNA copy number analysis results}
\label{fig:copy}
\end{figure}

We consider a real data example based on the DNA copy number analysis in \citet{hutter2007}.  In these applications, it is of biological importance to identify the change points, so the proposed method would be useful.  Data on the copy number for a particular gene are displayed in grey dots in Figure~\ref{fig:copy}(a).  We fit the proposed empirical Bayes model to these data, using the plug-in estimator for the model variance, which in this case is $\hat\sigma^2 = 0.093$, just like in Table~2 of \citet{hutter2007}.  Plot of the posterior mean estimate and is also shown.  The fit here appears to be quite good, perhaps with the exception around 600, and arguably the reason for this is the within-group variance seems to be much larger here than in other regions.  Interestingly, the distribution of $|B|$ in Figure~\ref{fig:copy}(b) is concentrated on much smaller values than in \citet{hutter2007}, who estimates about 15 piecewise constant blocks.  But a simple visual inspection of the data suggests much fewer blocks, and roughly 6--7 seems much more reasonable than 15.  

\subsection{Eye movement signal analysis}
Another interesting application of our method when $K=1$ is eye movement signal denoising. Eye movement of human and other foveate animals when scanning scenes is characterized by a fixate-saccade-fixate pattern. During the fixation phase, gaze position stables on the order of 0.2-0.3 seconds; in the saccadic phase, eye moves quickly on the order of 0.01-0.1 seconds. The time series of gaze position in terms of vertical and horizontal visual angle degree can be well approximated by piecewise linear functions, under the assumption that eye moves at an approximately constant velocity during either fixation or saccade phase; see \citet{pekkanen2017new}. 

\begin{figure}[t]
\centering
    \begin{subfigure}[b]{0.5\textwidth}
                \centering
                \includegraphics[width=0.99\textwidth]{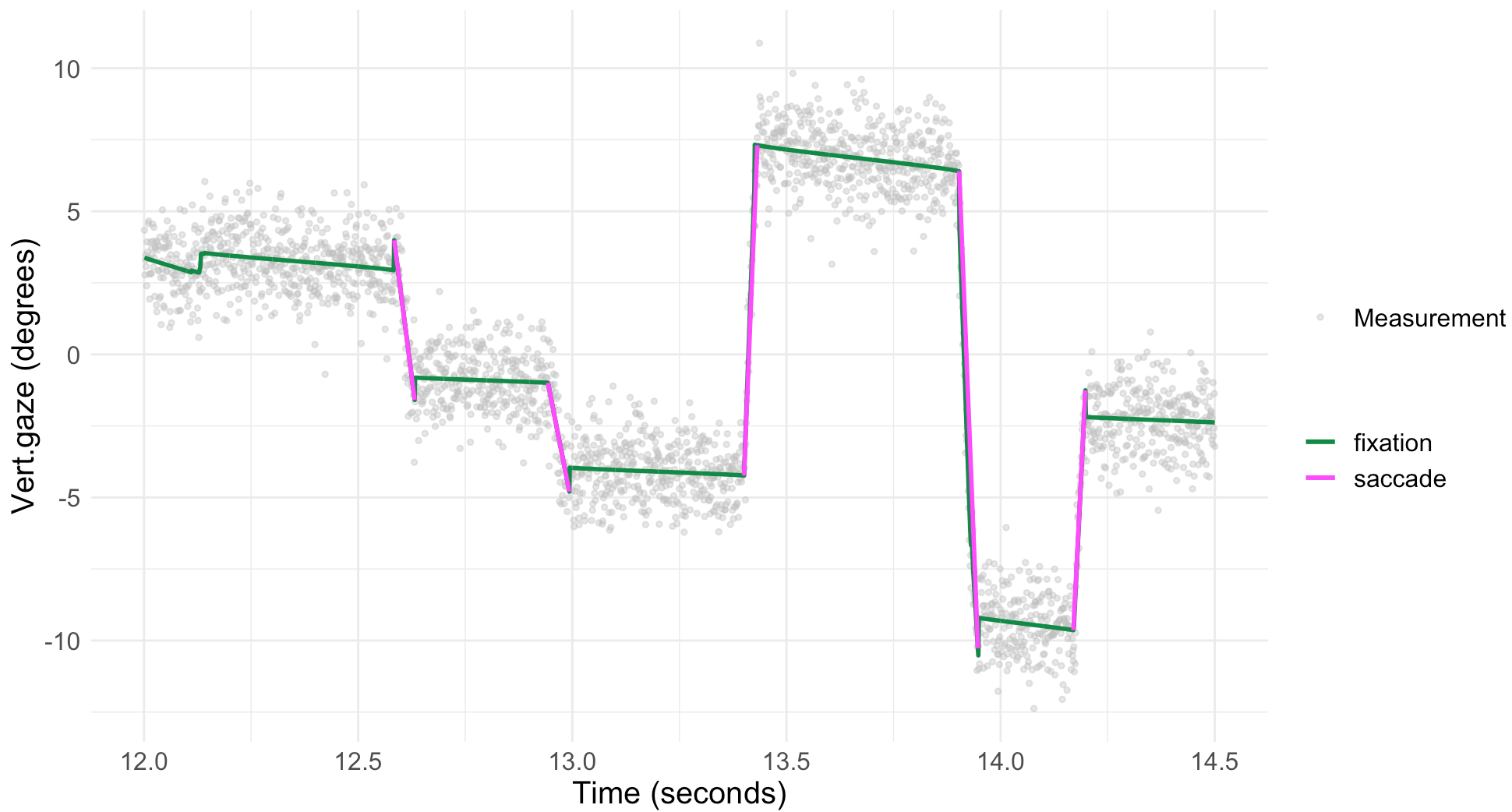}
        \caption{Vertical gaze position}
    \end{subfigure}%
    \begin{subfigure}[b]{0.5\textwidth}
                    \centering
                \includegraphics[width=0.99\textwidth]{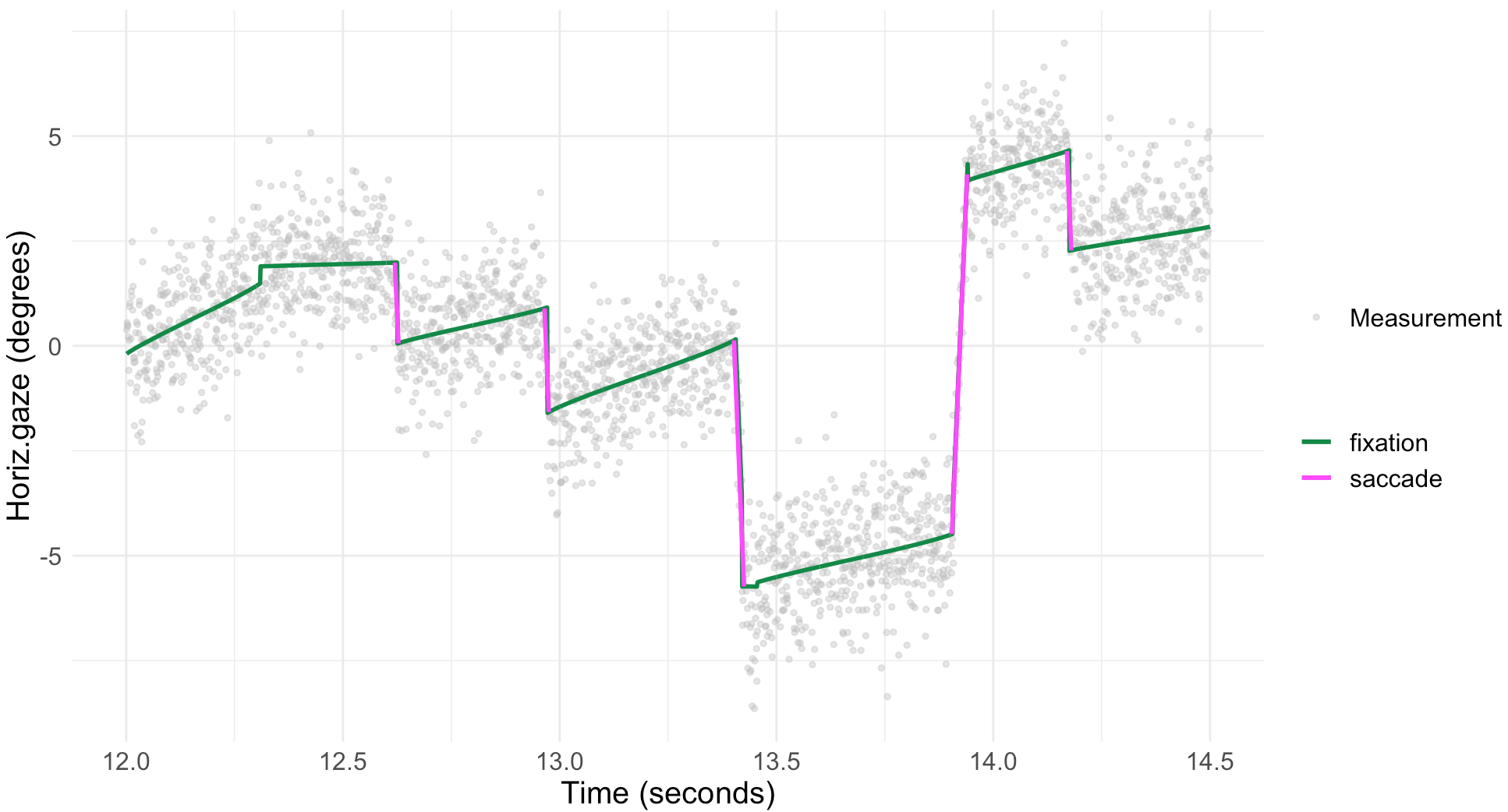}
         \caption{Horizontal gaze position}
    \end{subfigure}%
\caption{Eye movement signal denoising results}
\label{fig:eye}
\end{figure}

\begin{figure}[t]
\centering
    \includegraphics[width=0.60\textwidth]{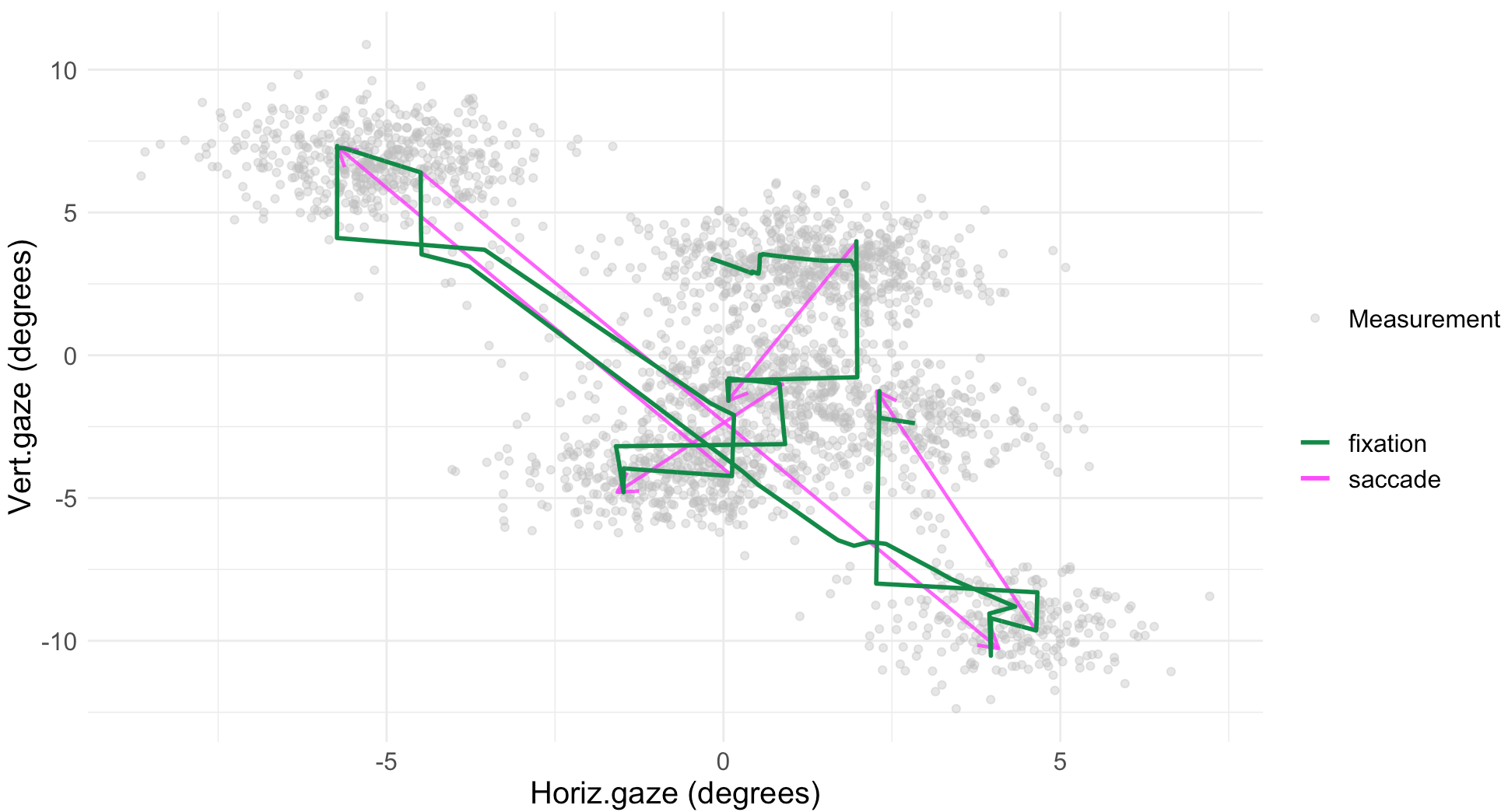}
                \caption{Estimated gaze path}
\label{fig:path}
\end{figure}

Noise in eye-movement recording is usually inevitable, ranging from around $0.01^\circ$ with laboratory optical equipment to well over $1^\circ$ in mobile recording with moving cameras. Here we consider the gaze position dataset in \citet{vig2012space} when participants watch a movie clip. Recording noise level is not reported in \citet{vig2012space}, so we adopt the  procedure in \citet{pekkanen2017new} where the same dataset is investigated, by adding a simulated measurement noise with standard deviation $1^\circ$; see Supplementary Information in \citet{pekkanen2017new}. Then for both vertical and horizontal gaze position data in a 2.5-second excerpt of the full recording, we fit an empirical Bayes estimator for the mean gaze trajectory with $\lambda=1$ and 50000-length MCMC after burn-in. The posterior mean estimate and the measurements mimicking mobile recording using a moving camera are plotted in Figure~\ref{fig:eye}. Based on the fitted vertical gaze position and horizontal gaze position, an estimated mean gaze path is plotted in Figure~\ref{fig:path}.

Our method helps us to identify and understand the segmentation of  fixate-saccade-fixate pattern in eye movements. As can be seen from Figures~\ref{fig:eye}--\ref{fig:path}, in the fixation segments (coloured in dark green), the eye moves slowly and steadily, hence the gaze position appears to be linear with a slope close to $0$; in the saccade segements (coloured in magenta), the gaze position is still linear but much steeper, showing a jump pattern. In addition, segmentation of eye movements is consistent between vertical gaze signal and horizontal gaze signal.

\section{Concluding remarks}
\label{S:discuss}

This paper considered inference on a piecewise polynomial signal where the degree is known but the block structure is unknown.  We developed an empirical Bayes posterior distribution that is simple and fast to compute and accompanied by a range of desirable theoretical results, including optimal posterior concentration rates and block selection consistency.
The general results are new and, when specialized to cases that have been investigated by others in the literature, our assumptions and/or conclusions are as good or better than what is currently available.  And, as our numerical results demonstrate, the strong theoretical properties of the proposed method carry over to real applications, where we see a considerable improvements compared to trend filtering, in particular, in cases where the underlying function being estimated is discontinuous or approximately so, like in Model~\ref{model:model6} above.  

\begin{figure}[t]
\begin{center}
\scalebox{0.55}{\includegraphics{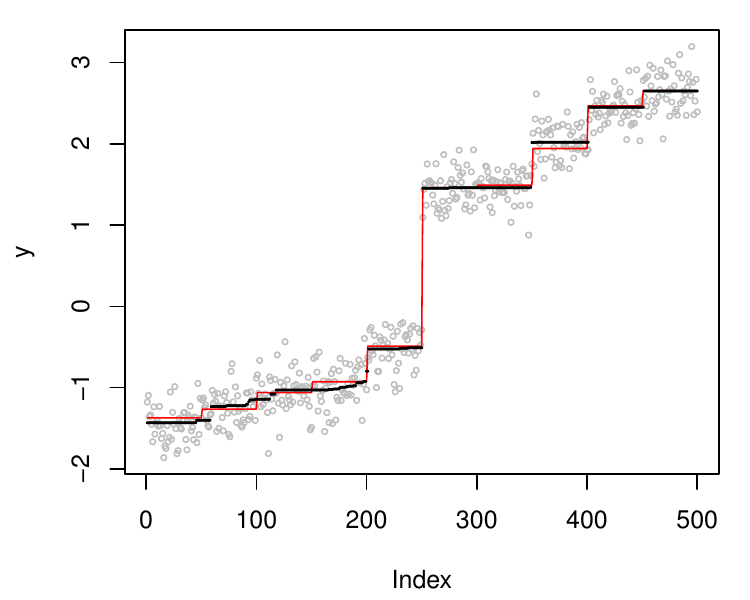}}
\end{center}
\caption{Simulated data (gray) from a distribution with a monotone piecewise constant mean (red) and the projection posterior mean (black).}
\label{fig:mono}
\end{figure}

There has been recent interest in cases where the signal is both piecewise constant and {\em monotone}; see, e.g., \citet{gao2017minimax} and \citet{guntuboyina.sen.2018}.  Of course, the method developed above can be applied in cases where the piecewise constant signal is monotone, but it is not immediately clear how to incorporate monotonicity into the prior formulation directly.  A clever alternative strategy to force the monotonicity constraint by projecting the posterior samples of $\theta$ from $\Pi^n$ onto the space of monotone sequences.  That is, if $\theta \sim \Pi^n$, then set 
\[ \text{proj}(\theta) = \arg\min_{z \in \Theta^\uparrow} \|z-\theta\|, \]
where $\Theta^\uparrow \subset \RR^n$ is the set of monotone sequences.  This projection operation is just a function of $\theta$, albeit implicit, so there is a corresponding posterior distribution for the projection, which is called a {\em projection posterior}.  General details about the projection posterior can be found in \citet{chakraborty.ghosal.proj}.  Aside from inheriting many of the desirable properties of the original posterior, the projection posterior is also relatively simple to compute.  The R package ``Iso'' \citep{R.iso} contains an implementation of the ``pool adjacent violators algorithm,'' or {\em PAVA} for short.  So, all we have to do is generate samples of the piecewise constant $\theta$ from the posterior $\Pi^n$ and then apply the {\tt pava} function to project it onto the space of monotone sequences.  Figure~\ref{fig:mono} shows the results of sampling from this projected posterior for a simulated data set, and the corresponding estimate appears to be quite accurate.  

Another interesting possible extension of the work here is related to the formulation in \citet{fan2018approximate}.  Consider a graph $G=(V,E)$ and, at each vertex $i \in V$, there is a response $Y_i \sim \nm(\theta_i^\star, \sigma^2)$, but only a small number of edges $(i,j) \in E$ have $\theta_i^\star \neq \theta_j^\star$.  That paper derives bounds on the recovery rate analogous to those achieved here in the chain graph/sequence model.  The only obstacle preventing us from extending our analysis to this more general setting is the need to assign a prior distribution for the block structure $B$ in this more complex graph.  For example, in a two-dimensional lattice graph, as might be used in imaging applications, one would need a prior on all possible ways that the lattice can be carved up into connected chunks, which seems non-trivial.  But given such a prior, we expect that the theoretical results described here would carry over directly.

\section*{Acknowledgments}

This work is partially supported by the U.S.~National Science Foundation, grants DMS--1737929, DMS--1737933, and DMS--1811802, and by the Simons Foundation Award 512620. The authors thank Marcus Hutter for sharing the
DNA copy number data.

\appendix

\section{Proofs}
\label{S:proofs}

\subsection{Preliminary results}

Before getting to the proofs of the main theorems, we first present some preliminary lemmas that helps to construct the proofs. First, we define 
\[ A_n=\{\theta \in \Theta_n(K): \|\theta-\theta^\star\|^2_n > M_n \eps_n(\theta^\star)\}.\]
For our theoretical analysis, it will help to rewrite the posterior distribution $\Pi^n(A_n)$ as the ratio $\Pi^n(A_n) = N_n(A_n) / D_n$. The numerator and denominator are 
\begin{align*}
N_n(A_n) & = \sum_B \pi_n(B) \int_{A \cap \Theta^B} R_n(\theta^B)^\alpha \, \Pi_n(d \theta^B \mid B),  \\
D_n & = \sum_B \pi_n(B) \int_{\Theta^B} R_n(\theta^B)^\alpha \, \Pi_n(d\theta^B \mid B), 
\end{align*}
where $R_n(\theta) = L_n(\theta)/L_n(\theta^\star)$ is the likelihood ratio, $\Theta_B$ is the column space of $Z^B$, i.e.,  $\Theta_B$
consists of all $n$-dimensional vectors $\theta^B$ such that
$\theta^B=Z^B\beta^B$. For a properly chosen matrix $Z$, $N_n(A_n)$ and $D_n$ can be rewritten in terms of $\beta^B$, given $B$,
\begin{align}
N_n(A_n) & = \sum_B \pi_n(B) \int_{\{\beta^B \in \RR^{|B|(K+1)}: Z^B\beta^B \in A_n\}} R_n(Z^B\beta^B)^\alpha \, \tilde{\pi}_n(\beta^B\mid B)\, d \beta^B, \label{eq:numerator}
\\
D_n & = \sum_B \pi_n(B) \int_{\RR^{|B|(K+1)}} R_n(Z^B\beta^B)^\alpha \, \tilde{\pi}_n(\beta^B\mid B) \,d \beta^B. \label{eq:denominator}
\end{align}
For the given $Z$ and $B_{\theta^\star}$, we let $\beta^\star$ be such that $\theta^\star=Z^{B_{\theta^\star}}\beta^\star$ and let $\beta_s^\star$ be its $s^\text{th}$ $(K+1)$-dimensional component.   We also abbreviate $B_{\theta^\star}$ by $B^\star$ and $\eps_n^2(\theta^\star)$ by $\eps_n^2$.  As discussed above, $B^\star$ may not be unique, but certain features of $B^\star$ are determined, in particular, its size $|B^\star|$ which, in turn, determines other features like $\pi_n(B^\star)$, etc.

\begin{lem}
\label{lem:den}
There exists a constant $c=c(\alpha, \sigma^2, K)$ such that $D_n \gtrsim \pi_n(B^\star)e^{-c|B^\star|}$ for all sufficiently large $n$.
\end{lem}

\begin{proof}
Given that $D_n$ is a sum of non-negative terms, it is straightforward to have
\[D_n > \pi_n(B^\star) \int_{\RR^{|B^\star|(K+1)}} R_n(Z^{B^\star}\beta^{B^\star})^\alpha \, \tilde{\pi}_n(\beta^{B^\star}\mid B^\star)\, d \beta^{B^\star}. \]
The integral in the right-hand side of the above inequality can be further written as,
\begin{equation*}
\prod_{s=1}^{|\Bs|}\int
   e^{-\frac{\alpha}{2\sigma^2}\{\|\bY_{\Bs(s)}-\ZBs\bbeta^{B^\star}_{s}\|^2-\|\bY_{\Bs(s)}-\ZBs\bbeta^\star_{s}\|^2\}} \N(\bbeta^{B^\star}_{s} \mid \hat{\bbeta}^{B^\star}_{s}, v(\ZBs^\top\ZBs)^{-1}) \, d \bbeta^{B^\star}_{s}.
\end{equation*}
Direct calculation shows that the above quantity equals
\[e^{\frac{\alpha}{2\sigma^2}\|\ZBs(\hat{\bbeta}_{s}^{B^\star}-\bbeta_{s}^\star)\|^2}\big(1+\tfrac{\alpha v}{\sigma^2}\big)^{-\frac{(K+1)|\Bs|}{2}} \geq \big(1+\tfrac{\alpha v}{\sigma^2}\big)^{-\frac{(K+1)|\Bs|}{2}} .\]
Therefore, $D_n > \pi_n(\Bs)e^{-c|\Bs|}$, where $c=\tfrac{(k+1)}{2}\log(1+\frac{\alpha v}{\sigma^2}) >0$. 
\end{proof}

\begin{lem}
\label{lem:num}
Take $q>1$ such that $\alpha q <1$. Then $\E_{\theta^*}\{N_n(A_n)\} \lesssim e^{-M_nr n\eps_n^2}$, for all large $n$, where $r=\alpha(1-q\alpha)/2\sigma^2$.
\end{lem}

\begin{proof}
Towards an upper bound, we interchange expectation with the finite sum over $B$ and the integral over $\beta^B$, the latter step justified by Tonelli's theorem, so that 
\begin{equation}
\label{eq:num.expectation}
\E_{\theta^\star}\{N_n(A_n)\} = \sum_B \pi_n(B) \int_{\{\beta^B: Z^B\beta^B \in A_n\}} \E_{\theta^\star}\{ R_n(Z^B\beta^B)^\alpha \, \tilde{\pi}_n(\beta^B\mid B)\}\,d \beta^B. 
\end{equation}
Next, we work with each of the $B$-dependent integrands separately.  For $q > 1$ such that $\alpha q < 1$, define the H\"older conjugate $p=q/(q-1) > 1$.  Then H\"older's inequality gives 
\[ \E_{\theta^\star}\{ R_n(Z^B\beta^B)^\alpha \, \tilde{\pi}_n(\beta^B\mid B)\} \leq \E_{\theta^\star}^{1/q} \{R_n(Z^B\beta^B)^{\alpha q}\} \, \E_{\theta^\star}^{1/p}\{ \tilde{\pi}_n(\beta^B\mid B)^p \}. \]
On the set $A_n$, since $\alpha q < 1$, the first term above is uniformly bounded by $e^{-M_n rn \eps_n^2}$.  To see this, note that, for a general $\theta \in A_n$, if $p_\theta^n$ denotes the joint density of $Y$ under \eqref{eq:model}, and $D_{\alpha q}$ the R\'enyi $\alpha q$-divergence of one normal distribution from another \citep[e.g.,][p.~3800]{van2014renyi}, then  
\begin{align*}
\E_{\theta^\star}\{R_n(Z^B\beta^B)^{\alpha q}\} & = \int \{ p_{\theta^B}^n(y) \}^{\alpha q} \, \{ p_{\theta^\star}^n(y) \}^{1-\alpha q} \,dy 
= e^{-\frac{\alpha q(1-\alpha q)n}{2\sigma^2} \|\theta^B-\theta^\star\|_n^2}. 
\end{align*}
Then for the second term in the upper bound above, using prior in \eqref{eq:prior_beta} we can have
\begin{equation}
    \E_{\theta^\star}^{1/p}\{ \tilde{\pi}_n(\beta^B\mid B)^p\}=\prod_{s=1}^{|B|} \E^{1/p}_{ {\theta}^*}\{\N^p(\beta_s^B \mid\hat{ {\beta}}^B_{s}, v(Z_{B(s)}^\top Z_{B(s)})^{-1})\},
    \label{eq:dtheta}
\end{equation}
which is equivalent to
\begin{align}
    \prod_{s=1}^{|B|} 
    \frac{|Z_{B(s)}^\top Z_{B(s)}|^{\frac{1}{2}}}{(2\pi v)^{\frac{K}{2}}}
    \E^{1/p}_{ {\theta}^\star}(e^{-\frac{p\sigma^2}{2v}V_{s}}),
    \label{eq:chisq}
\end{align}
where 
\[ V_{s}=\sigma^{-2}\|P_{B(s)}(Z_{B(s)} {\beta}^B_{s}- {Y}_{B{(s)}})\|^2,\]
is distributed as a non-central chi-square with $(K+1)$ degrees of freedom and non-centrality parameter
\[ \lambda_{s}=\sigma^{-2}\|P_{B(s)}(Z_{B(s)} {\beta}_{s}^B- {\theta}^\star_{B(s)})\|^2, \]
with $P_{B(s)}=Z_{B(s)}(Z_{B(s)}^\top Z_{B(s)})^{-1}Z_{B(s)}^\top$.  Using the familiar formula for the moment generating function of non-central chi-square, we have
\begin{align}
     \E^{1/p}_{ {\theta}^*}(e^{-\frac{p\sigma^2}{2v}V_{s}})&=(1+p\sigma^2/v)^{-K/(2p)}e^{-\frac{1}{2(v/\sigma^2+p)}\lambda_s}
     \label{eq:mgf}
\end{align}
For the integral
\[\int_{\RR^{|B|(K+1)} } \E_{\theta^\star}^{1/p}\{ \tilde{\pi}_n(\beta^B\mid B)^p\} \, d \beta^B,\]
if we plug \eqref{eq:dtheta}, \eqref{eq:chisq}, and \eqref{eq:mgf} into the integrand, then it simplfies as
\begin{align*}
    \prod_{s=1}^{|B|} \tfrac{|Z_{B(s)}^\top Z_{B(s)}|^{\frac{1}{2}}}{(2\pi v)^{\frac{K}{2}}} & (1+\tfrac{p\sigma^2}{v})^{-\frac{K}{2p}} \\
    & \times 
    \int_{\RR^{K+1}} \N( {\beta}^B_{s}\mid(Z_{B(s)}^\top Z_{B(s)})^{-1}Z_{B(s)}^\top {\theta}^\star_{B(s)}, (\tfrac{v}{\sigma^2}+p)(Z_{B(s)}^\top Z_{B(s)})^{-1}) \,d \beta_s^B.
\end{align*}
By direct calculation, this can be written as $\zeta^{|B|}$, where $\zeta=\big\{\tfrac{(1+p\sigma^2/v)^{1/q}}{\sigma^2}\big\}^{(K+1)/2}$.  Therefore, we have 
$$\E_{\btheta^\star}\{N_n(A_n)\} \le e^{-M_nrn\eps_n^2}\sum_B\zeta^{|B|}\pi_n(B)=e^{-M_nrn\eps_n^2}\sum_{b=1}^n\zeta^bf_n(b).$$
Given that $f_n(b)\propto n^{-\lambda (b-1)}$ and $\zeta$ is a positive constant, the summation term in the above upper bound is uniformly bounded in $n$, proving the claim.
\end{proof}

\subsection{Proof of Theorem~\ref{thm:rate}}
\label{SS:proofs.rate}

By Lemma~\ref{lem:den}, for sufficiently large $n$, we have 
\[\E_{\btheta^\star}\{\Pi^n(A_n)\} \le \frac{e^{c|B^\star|}}{\pi_n(B^\star)}\E_{\btheta^\star}\{N_n(A_n)\}.\]
Plug in the bound from Lemma~\ref{lem:num} to get
\[\E_{\btheta^\star}\{\Pi^n(A_n)\} \lesssim e^{c|B^\star|-M_n r n\eps_n^2} \frac{\binom{n-1}{|B^\star|-1}}{f_n(B^\star)}.\]
On the one hand, if $|B_{\theta^\star}|=1$, then both $n\eps_n^2$ and the ratio in the above display are constant.  Therefore, the upper bound is $\lesssim e^{-G M_n} \to 0$ for a constant $G$ and $M_n \to \infty$.  On the other hand, if $|B_{\theta^\star}| \geq 2$, then $n\eps_n^2$ is diverging, and we take $M_n \equiv M$ a fixed constant.  Also, using the formula for $f_n(|B^\star|)$ in \eqref{eq:S.prior} and the standard bound, $\binom{n}{b} \leq e^{b \log(en/b)}$, on the binomial coefficient, we get
\begin{equation}
\label{eq:thm1.bound}
\E_{\theta^\star}\{\Pi^n(A_n)\} \lesssim \exp\{-M r n \eps_n^2 + n\eps_n^2 + \lambda |B^\star| \log n + c|B^\star|\}. 
\end{equation}
The exponent on the right-hand side can be rewritten as 
\[ -n \eps_n^2 \Bigl( Mr - 1 - \frac{\lambda |B^\star| \log n}{n \eps_n^2} - \frac{c|B^\star|}{n \eps_n^2} \Bigr). \]
Since $|B^\star| = o(n)$, if we let $n \eps_n^2 = |B^\star| \log n$, the two ratios inside the parentheses in the above display are upper bounded by a constant for sufficiently large $n$.  Therefore, for a sufficiently large $M$, there exists a constant $G > 0$, depending on the constants $M$, $r$, $\lambda$, and $c=c(\alpha, \sigma^2, K)$, such that the right-hand side of \eqref{eq:thm1.bound} can be written as $e^{-G n \eps_n^2} \to 0$.

\subsection{Proof of Theorem~\ref{thm:mean}}
\label{SS:proofs.mean}

It follows from Jensen's inequality that $\|\hat\theta - \theta^\star\|_n^2 \leq \int \|\theta - \theta^\star\|_n^2 \, \Pi^n(d\theta)$.  So it suffices to bound the expectation of the upper bound.  Towards this, write $\RR^n$ as $A \cup A^c$, where $A=A_n$ is as defined above.  Then 
\begin{equation}
\label{eq:mean.bound}
\E_{\theta^\star} \int \|\theta-\theta^\star\|_n^2 \, \Pi^n(d\theta) \leq M_n \eps_n^2 + \E_{\theta^\star} \int_A \|\theta-\theta^\star\|_n^2 \, \Pi^n(d\theta).
\end{equation}
That remaining integral can be expressed as a ratio of numerator to denominator, where the denominator $D_n$ is just as in Lemma~\ref{lem:den} and the numerator $\widetilde N_n(A)$ is 
\begin{align*}
\widetilde N_n(A) & = \int_A \|\theta-\theta^\star\|_n^2 \, R_n(\theta)^\alpha \, \Pi_n(d\theta) \\
& = \sum_B \pi_n(B) \int_{A \cap \Theta_B} \|\theta^B-\theta^\star\|_n^2 \, R_n(\theta^B)^\alpha \, \pi_n(d\theta^B \mid B). 
\end{align*}
Take expectation of the numerator to the inside of the integral and apply H\"older's inequality just like in the proof of Lemma~\ref{lem:num}.  This gives the following upper bound on each $B$-specific integral:
\[ \int_{A \cap \Theta_B} \|\theta-\theta^\star\|_n^2 e^{-h\|\theta-\theta^\star\|^2} \E_{\theta^\star}^{1/p} \{ \pi_n(d\theta^B \mid B)^p\}, \]
where $h > 0$ is a constant that depends only on $\alpha$, $\sigma^2$, and the H\"older constant $q > 1$.  Since the function $x \mapsto x e^{-h x}$ is eventually monotone decreasing, for sufficiently large $M$ we get a trivial upper bound on the above display, i.e., 
\[ M_n \eps_n^2 e^{-M_n h n\eps_n^2} \int_{\Theta_B} \E_{\theta^\star}^{1/p} \{ \pi_n(d\theta^B \mid B)^p\}. \]
The same argument as above bounds the remaining integral by $\zeta^{|B|}$, and the prior $\pi_n(B)$ takes care of that contribution.  In the case where $|B^\star|=1$, where $n\eps_n^2$ is bounded, choosing $M_n \to \infty$ will take care of the bound on $D_n$ from Lemma~\ref{lem:den}.  Similarly, for cases when $n\eps_n^2 \to \infty$, we can use a sufficiently large constant $M_n \equiv M$ to take care of the lower bound on $D_n$.  Therefore, the second term on the right-hand side of \eqref{eq:mean.bound} is also bounded by a multiple of $M_n \eps_n^2$.

\subsection{Proof of Theorem~\ref{thm:dim}}
\label{SS:proofs.dim}

Following arguments similar to those in the proof of Theorem~\ref{thm:rate}, we get 
\[\E_{\theta^\star} \pi_n(\{B: |B| > C|B^\star|\}) \lesssim \frac{e^{c|B^\star|}}{\pi_n(B^\star)}\sum_{b=C|B^\star|+1}^n \zeta^b f_n(b).
\]
From the formula \eqref{eq:S.prior} for $f_n$, factor out a common $(\zeta n^{-\lambda})^{C|B^\star|}$ from the summation, which will be the dominant term.  Indeed, like in the proof of Theorem~\ref{thm:rate}, the ratio in the above display is of order $\exp\{n\eps_n^2 + \lambda |B^\star|\log n\}$.  Then the right-hand side above is of order 
\[ \exp\{n\eps_n^2 + \lambda |B^\star|\log n + C|B^\star|\log \zeta - C\lambda |B^\star|\log n\}. \]
The exponent can be written as 
\[ -|B^\star| \log n \Bigl(C\lambda - \lambda - \frac{C\log \zeta}{\log n} - \frac{n\eps_n^2}{|B^\star|\log n} \Bigr). \]
Since $n\eps_n^2 \leq |B^\star|\log n$, it is clear that, if $C$ is strictly larger than $1 + \lambda^{-1}$, then the term in parentheses is bigger than some constant $G > 0$ for all large $n$.

\subsection{Proof of Theorem~\ref{thm:no.refinements}}

Choose and fix any $B^\star \in \BB^\star$. For a generic configuration $B$, we have 
\begin{align}
    \pi^n(B) & \le \pi^n(B) / \pi^n(B^\star) \notag \\
    & =\tfrac{\pi_n(B)}{\pi_n(\Bs)}(1+\tfrac{v\alpha}{\sigma^2})^{(K+1)(|\Bs|-|B|)/2}e^{-\frac{\alpha}{2\sigma^2}\{\sum_{s=1}^{|B|}\|(I-P_{B(s)})\bY_{B(s)}\|^2-\sum_{s=1}^{|B^\star|}\|(I-P_{\Bs(s)})\bY_{\Bs(s)}\|^2\}}\notag\\
    &=\tfrac{\pi_n(B)}{\pi_n(\Bs)}(1+\tfrac{v\alpha}{\sigma^2})^{(K+1)(|\Bs|-|B|)/2}e^{-\frac{\alpha}{2\sigma^2}\{\|(I_n-P^B)Y\|^2-\|(I_n-P^{\Bs})Y\|^2\}}
    \label{eq:post_marginal_bound}
\end{align}
where $P^B=Z^B(Z^{B\top} Z^B)^{-1}Z^{B\top}$, $P^{\Bs}=Z^{\Bs}(Z^{\Bs \top} Z^{\Bs})^{-1}Z^{\Bs \top}$ and $Z^B$ is defined in \eqref{eq:Z^B}.
If $B$ is a refinement of $B^\star$, then column space of $Z^{B^\star}$ is a subset of the column space of $Z^{B}$, i.e., $\mathcal{C}(Z^{B^\star}) \subseteq \mathcal{C}(Z^{B})$ and therefore, $P^B-P^{\Bs}$ is idempotent of rank $(K+1)(|B|-|\Bs|)>0$. Thus, \eqref{eq:post_marginal_bound} can be rewritten as,
\begin{equation}
    \pi^n(B) \le \frac{\pi_n(B)}{\pi_n(\Bs)}\Big(1+\frac{v\alpha}{\sigma^2}\Big)^{(K+1)(|B|-|\Bs|)/2}e^{\frac{\alpha}{2\sigma^2}V},
    \label{eq:central_chisq}
\end{equation}
where $V=\bY^\top(P^B-P^{\Bs})\bY$ is distributed as a central chi-square with $(K+1)(|B|-|B^\star|)$ degrees of freedom. From the chi-square moment generating function, we get
\[\E_{\btheta^\star}\{\pi^n(B)\} \le \frac{\pi_n(B)}{\pi_n(\Bs)}\psi^{|B|-|B^\star|},\]
where $\psi=\psi(\alpha, v, \sigma, K)$ is a positive constant. For a suitable constant $C$ as in Theorem~\ref{thm:dim}, write $\B_n = \{B: B \sqsupset B^\star, \, |B| \leq C|B^\star|\}$.  Then
\[ \{B: B \sqsupset B^\star\} \subseteq \B_n \cup \{B: |B| > C|B^\star|\}. \]
Since the right-most event has vanishing $\Pi^n$-probability by Theorem~\ref{thm:dim}, it follows that we can focus just on the event $\B_n$, and   
\[ \E_{\theta^\star} \pi^n(\B_n) = \sum_{B \in \B_n} \E_{\theta^\star} \pi^n(B)  \leq \sum_{B \in \B_n} \frac{\pi_n(B)}{\pi_n(B^\star)} \psi^{|B|-|B^\star|}. \]
Plug in the prior for $B$---which only depends on $|B|$---and simplify:
\[  \E_{\theta^\star} \{ \Pi^n(\theta: B_\theta \in \B_n)\} \leq \sum_{b = b^\star + 1}^{C b^\star} \frac{\binom{n-1}{b^\star-1} \binom{n-b^\star}{b-b^\star}}{\binom{n-1}{b-1}} (\psi n^{-\lambda})^{b-b^\star}, \]
where $b^\star = |B^\star|$.  From 
\[ \frac{\binom{n-1}{b^\star-1} \binom{n-b^\star}{b-b^\star}}{\binom{n-1}{b-1}} = \binom{b-1}{b^\star-1} \leq b^{b-b^\star}, \]
and the assumption that $b^\star = o(n^\lambda)$, we get that the summation on the right-hand side above is upper-bounded by 
\[ \sum_{b = b^\star+1}^{Cb^\star} (\psi b n^{-\lambda})^{b-b^\star} \leq \sum_{b=b^\star + 1}^{Cb^\star} (C \psi b^\star n^{-\lambda})^{b-b^\star} \lesssim e^{-\lambda b^\star \log n}, \quad \text{for all large $n$}. \]
This argument can be duplicated for any $B^\star \in \BB^\star$, and there are at most $O(|B^\star|)$ many equivalent block configurations, so we get 
\[ \E_{\theta^\star} \pi^n(\{B: B \sqsupset B^\star \text{ for some $B^\star \in \BB^\star$}\}) \lesssim |B^\star| e^{-\lambda |B^\star| \log n} \to 0. \]

\subsection{Proof of Theorem~\ref{thm:B_consist}}

Choose and fix any $B^\star \in \BB^\star$.  In light of Theorem~\ref{thm:no.refinements}, it suffices to show that
\[\E_{\theta^\star}\pi^n(\{B: B \not\sqsupseteq B^\star\}) \to 0.\] 
For a generic $B \not\sqsupseteq B^\star$, according to \eqref{eq:post_marginal_bound} we have
\[\pi^n(B) \le \frac{\pi_n(B)}{\pi_n(\Bs)}\Big(1+\frac{v\alpha}{\sigma^2}\Big)^{(K+1)(|B|-|\Bs|)/2}e^{\frac{\alpha}{2\sigma^2}Y^\top(P^B-P^{B^\star})Y}.\]
We proceed with a proof for the piecewise constant ($K=0$) case first, then describe how the general $K$ case is the same.  Let $\theta$ be a piecewise constant signal corresponding to the block configuration $B$. 
Then we can rewrite $\theta$ as $X\eta$, where $X$ is an $n\times n$ lower triangular matrix with unit entries and 
\[\eta=(\theta_1, \theta_2-\theta_1,\theta_3-\theta_2,\ldots,\theta_n-\theta_{n-1})^\top.\]
It is easy to show that $\eta$ is sparse. Let $J=\{1\}\cup\{j:\eta_j\neq 0\}$, then $|J|=|B|$. Let $\eta_J$ be the $|J|$-vector containing the particular  entries with their indices in $J$ and $X_J$ be the columns of $X$ corresponding to $J$. Then we can also write $\theta=X_J\eta_J$. Hence we can reformulate model \eqref{eq:model} as 
\begin{equation}
    Y=X\eta^\star+\sigma \xi, \quad \xi \sim \N(0, I).
    \label{eq:regression}
\end{equation}
Under this formulation, recovering block structure $B^\star$ is equivalent to identifying the non-zero coefficients in $\eta^\star$, i.e., recovering $J^\star$. One basic observation is that $P^B$ is equal to $P_J=X_J (X_J^\top X_J)^{-1}X_J^{-1}$. Then we can rewrite  $Y^\top(P^B-P^{B^\star})Y$ as, 
\[-\|(I-P_J)X\eta^\star\|^2-2\sigma \xi^\top(I-P_J)X\eta^\star+\sigma^2\xi^\top(P_J-P_{J^\star})\xi.\]
In addition, because $(P_{J^\star}-P_{J \cap J^\star})$ is positive definite, 
the right-most quadratic form above can be bounded as follows,
\[\xi^\top(P_J-P_{J^\star})\xi=\xi^\top(P_J-P_{J \cap J^\star})\xi-\xi^\top(P_{J^\star}-P_{J \cap J^\star})\xi
\le \xi^\top(P_J-P_{J \cap J^\star})\xi.\]
 Therefore, $Y^\top(P^B-P^{B^\star})Y$ can be bounded above by,
\[-\|(I-P_J)X\eta^\star\|^2-2\sigma \xi^\top(I-P_J)X\eta^\star+\xi^\top(P_J-P_{J \cap J^\star})\xi.\]
Note that the second and third terms in the above upper bound follow normal and chi-square distributions respectively, and additionally $(I-P_J)(P_J-P_{J \cap J^\star})=0$ implies independence.
Hence, using normal and chi-square moment generating functions we can have,
\[\E_{\theta^\star}[e^{\frac{\alpha}{2\sigma^2}Y^\top(P^B-P^{B^\star})Y}] \le (1-\alpha)^{-\frac{1}{2}(|J^\star|-|J \cap J^\star|)}e^{-\frac{\alpha(1-\alpha)}{2\sigma^2}\|(I-P_J)X\eta^\star\|^2}.\]
Since 
\[\|(I-P_J)X\eta^\star\|^2=\|(I-P_J)X_{J^\star\cap J^c} \eta^\star_{J^\star\cap J^c}\|,\]
if we let
\[\delta_n=\min_{j \in J^\star\cap\{j>1\}}|\theta^\star_j-\theta^\star_{j-1}|,\]
then
\begin{align*}
\|(I-P_J)X\eta^\star\|^2 & \ge \lambda_{\min}\big(X_{J^\star\cap J^c}^\top X_{J^\star\cap J^c}\big)\|\eta^\star_{J^\star\cap J^c}\|^2 \\
& \ge \lambda_{\min}\big(X_{J^\star}^\top X_{J^\star}\big) \delta_n^2 (|J^\star|-|J^{\star}\cap J|).
\end{align*}
Then we let 
\[\gamma_n=\min_{j,j'\in J^\star, j \neq j'} |j-j'|,\]
according to Lemma~\ref{lem:min.B} below,
$\|(I-P_J)X\eta^\star\|^2$ can be further lower bounded by,
\[\gamma_n\delta_n^2(|J^\star|-|J^\star \cap J|)/4.\]
Therefore, if we let
\[\gamma_n\delta^2_n\ge \frac{4M\sigma^2}{\alpha(1-\alpha)}\log n,\]
then 
\begin{equation}
    \E_{\theta^\star} \{\pi^n(B)\} \le \frac{\pi_n(B)}{\pi_n(\Bs)}\phi^{|\Bs|-|B|}(\omega n^{-M})^{|J^\star|-|J \cap J^\star|},
    \label{eq:pi_B_bound}
\end{equation}
where $\omega=(1-\alpha)^{-1/2}$ and $\phi=(1+\frac{v\alpha}{\sigma^2})^{-1/2}$.
Plug in the expressions for $\pi_n(B)$ and $\pi_n(\Bs)$ from \eqref{eq:S.prior} and then sum over all $B$ such that $B \not\sqsupseteq B^\star$ to get 
\[\E_{\theta^\star} \pi^n(\{B: B \not\sqsupseteq B^\star)\}) \le \sum_{b=1}^{C b^\star}\sum_{t=1}^{b \wedge b^\star} \frac{\binom{b^\star}{t}\binom{n-b^\star}{b-t} \binom{n}{b^\star}}{\binom{n}{b}}(\phi n^\lambda)^{b^\star-b}(\omega n^{-M})^{b^\star-t},\]
where $b^\star = |B^\star|$ and the first sum is restricted to $b \leq Cb^\star$ by Theorem~\ref{thm:dim}.  The ratio of binomial coefficients can be bounded as
\[\frac{\binom{b^\star}{t} \binom{n-b^\star}{b-t} \binom{n}{b^\star}}{\binom{n}{b}}=\binom{b}{t} \binom{n-b}{b^\star-t} \le (n^3)^{b \vee b^\star - t}.\]
Plug in this bound and split the sum over $b$ into two cases: $b \le b^\star$ and $b>b^\star$. For the first case, we have 
\[\sum_{b=1}^{b^\star}\sum_{t=1}^b(\phi n^\lambda)^{b^\star-b}(n^3)^{b^\star-t}(\omega n^{-M})^{b^\star-t}
\lesssim \sum_{b=1}^{b^\star}(\phi n^{3+\lambda-M})^{b^\star-b},\]
and the right-hand side vanishes since $M > 3 + \lambda$.  Similarly, for the second case
\[\sum_{b=b^\star+1}^{Cb^\star}\sum_{t=1}^{b^\star}(\phi n^\lambda)^{b^\star-b}(n^3)^{b^\star-t}(\omega n^{-M})^{b^\star-t}
\lesssim
\omega n^{3-M} \sum_{b=b^\star+1}^{Cb^\star} (n^{3-\lambda})^{b-b^\star},\]
and if $\lambda \ge 3$, then the sum is dominated by term $n^{3-M}$. In either case, the upper bound vanishes---actually the upper bound is $O(n^{-1})$ because $M > 4 + \lambda$---which proves the claim for the particular $B^\star \in \BB^\star$.  The above argument is not specific to any $B^\star$, so if we repeat the above argument sum over all such $B^\star \in \BB^\star$, then we get
\[ \sum_{B^\star \in \BB^\star} \E_{\theta^\star}\pi^n(\{B: B \not\sqsupseteq B^\star\}) \lesssim |B^\star| n^{-1} \to 0. \]
Finally, since 
\[ 1 - \E_{\theta^\star} \pi^n(\BB^\star) = \sum_{B^\star \in \BB^\star}  \E_{\theta^\star}\pi^n(\{B: B \sqsupset B^\star\}) + \sum_{B^\star \in \BB^\star} \E_{\theta^\star}\pi^n(\{B: B \not\sqsupseteq B^\star\}), \]
and first term on the right-hand side vanishes by Theorem~\ref{thm:no.refinements} and the second term vanishes by the argument above, we conclude that 
\[ \E_{\theta^\star} \pi^n(\BB^\star) \to 1, \quad n \to \infty. \]

{
Next, we show that for general piecewise polynomial $K \ge 1$, $\E_{\theta^\star}[\Pi^n(\theta:B_{\theta} \not\sqsupseteq B^\star)]$ can be bounded in a similar fashion. 
We define  
\[X^{(K)}=\begin{pmatrix}
I_{K} &  \\
    & L_{n-K}
\end{pmatrix}, \quad K=1,\ldots,n-1,\]
where $L_{n-K}$ is an $(n-K) \times (n-K)$-dimensional lower triangular matrix with unit entries.

Now, let's consider a generic degree-$K$ piecewise polynomial signal $\theta$ with underlying block configuration $B$, then $\theta$ can be written as
\[X \eta,\]
where $X=L_n X^{(1)}\cdots X^{(K)}$ and, 
\[\eta=\big((\Delta^0 \theta)_1, (\Delta ^1 \theta)_1,\ldots,(\Delta^K \theta)_{1}, \Delta^{K+1} \theta \big),\]
with $\Delta^{K}$ being the $K^\text{th}$-order difference operator defined in Section~\ref{SS:structure}. Note that $\eta$ is also sparse here, and if we let $J=\{1,\ldots,K+1\}\cup\{j:\eta_j \neq 0\}$, then $|J|=(K+1)|B|$. Therefore, a similar result to \eqref{eq:pi_B_bound} can be obtained,
\[\pi(B) \le \frac{\pi_n(B)}{\pi_n(\Bs)}\phi^{(K+1)(|B|-|\Bs|)/2}(\omega n^{-M})^{|J^\star|-|J \cap J^\star|}.\]
Then based on Lemma~\ref{lem:min.B}, using recursion, rest of the proofs can follow similar arguments in the $K=0$ case with 
\[\delta_n=\underset{j \in J^\star \cap \{j>K+1\}}{\min}|\eta_j|,\]
and 
\[\gamma_n\delta^2_n\ge \frac{4^{K+1}M\sigma^2}{\alpha(1-\alpha)}\log n.\]
}

\subsection{An eigenvalue bound}
\label{SS:min.B}

\begin{lem}
\label{lem:min.B}
Consider $J=\{j_1,\ldots,j_s\} \subset \{1,\ldots,n\}$, let $X$ be an $n \times n$-dimensional lower triangular matrix with unit entries and $X_S$ be the $n \times s$-dimensional sub-matrix of $X$ with the columns corresponding to $J$, define
\[\gamma= \min_{1 \leq \ell < k \leq s} |j_\ell-j_k|,\]
then the smallest eigenvalue of $X_S^\top X_S$ satisfies
\[\lambda_{\min}(X_S^\top X_S) > \gamma/4.\]
\end{lem}

\begin{proof}
Without loss of generality, here and throughout we assume that $j_1 <\cdots<j_s$. It is straightforward to observe that, 
\[X_S^\top X_S=
\begin{pmatrix}
a_1 & a_2 &\vdots& a_s\\
a_2 & a_2 & \vdots & \vdots\\
\cdots & \cdots& \cdot &\vdots\\
a_s & \cdots &\cdots &a_s
\end{pmatrix},\]
with $a_i=n+1-j_i$, $i=1,\ldots,s$. According to Lemma~3 in \citet{qian2016stepwise}, the inverse of $X^\top_S X_S$ is tridiagonal, i.e.,
\[(X_S^\top X_S)^{-1}=
\begin{pmatrix}
r_{11} & r_{12} &&&&\\
r_{21} & r_{22} & r_{23} &&&\\
&r_{32} & r_{33} & r_{34} &&\\
& & \ddots &\ddots &\ddots &\\
& & &r_{s-1,s-2} & r_{s-1,s-1} & r_{s-1,s}\\
&&&&r_{s,s-1}&r_{s,s}
\end{pmatrix}\]
where
\[r_{ij}=\begin{cases}
\frac{1}{a_1-a_2} & i=j=1\\
-\frac{1}{a_{j-1}-a_{j}}& i=j-1\\
-\frac{1}{a_j-a_{j+1}} &i=j+1\\
\frac{a_{j-1}-a_{j+1}}{(a_{j-1}-a_j)(a_j-a_{j+1})} & 1<i=j<s\\
\frac{a_{s-1}}{(a_{s-1}-a_s)(a_s)} & i=j=s\\
0 & \text{otherwise}.
\end{cases}
\]
For any vector $u \in \RR^{s\times 1}$,
\begin{align*}
    u^\top (X_S^\top X_S)^{-1} u & =
    \sum_{i=1}^s r_{ii} u_i^2 +
    \sum_{i=1}^{s-1} r_{i,i+1} u_i u_{i+1}+
    \sum_{i=1}^{s-1} r_{i+1,i} u_i u_{i+1}\\
    & \le \max_i |r_{ii}| \sum_{i=1}^s u_i^2+(\max_i|r_{i,i+1}|+ \max_i |r_{i+1,i}|) \sum_{i=1}^{s-1}|u_i u_{i+1}|
\end{align*}
Because $\sum_{i=1}^{s-1}|u_i u_{i+1}| \le \frac{1}{2}\sum_{i=1}^{s-1}(u_i^2+u_{i+1}^2) \le \sum_{i=1}^s u_i^2$, 
\begin{align*}
    u^\top (X_S^\top X_S)^{-1} u &\le \max_i |r_{ii}|+\max_i |r_{i,i+1}| + \max_i |r_{i+1,i}|)\sum_{i=1}^s u_i^2\\
    & \le \frac{4}{\gamma}\sum_{i=1}^s u_i^2.
\end{align*}
Thus, $\lambda_{\max}\{(X^\top_S X^\top_S)^{-1}\} \le 4/\gamma$, and therefore, $\lambda_{\min}(X^\top_S X^\top_S) \ge \gamma/4$.
\end{proof}

\bibliographystyle{apalike}
\bibliography{ref}

\end{document}